\input amstex
\documentstyle{amsppt}
\document


\magnification=1100 \NoBlackBoxes

\hsize=15cm \vsize=20cm \baselineskip=12pt

\pretolerance = 1600
\tolerance = 1600


\def\sv{{v}}        
\def\sw{{w}}        
\def\v{\nu}


\def\Chat{{\widehat{C}}}

\let\bb=\Bbb
\let\fr=\frak

\def\adots{\mathinner{\mkern1mu\raise1pt\vbox{\kern7pt\hbox{.}}\mkern2mu\raise4pt\hbox{.}\mkern2mu\raise7pt\hbox{.}\mkern1mu}}

\def\disp#1{\hskip40pt\hbox{(#1)}}


\def\bmatrix#1{\pmatrix #1 \endpmatrix}  


\def\\{\backslash}    
\def\A {{\bb A}}
\def\alf{\alpha}
\def\bra{\langle}
\def\<{\langle}

\def\C {{\bb C}}

\def\ch{\hbox{\rm ch}}

\def\eps{\varepsilon}

\def\F {{\bb F}}

\def\gam{\gamma}

\def\isom{\approx}

\def\J{{\bb J}}

\def\K {{\bb K}}

\def\ket{\rangle}
\def\>{\rangle}

\def\meas{\hbox{\rm meas}\,}

\def\N {{\cal N}}

\def\o {{\fr o}}

\def\ph{\varphi}

\def\Q {{\bb Q}}

\def\R {{\bb R}}

\def\th{^{\scriptstyle th}}

\def\to{\longrightarrow}

\def\Z {{\bb Z}}

\def\tr{{\text tr}}



\def\pht{\widetilde{\ph}}

\def\Pe{\hbox{P\'e}}

\def\diff{{\fr d}}

\def\hf{{\scriptstyle {1\over 2}}}

\def\Hf{{\textstyle{1\over 2}}}

\def\chibar{{\overline{\chi}}}

\def\psibar{\overline{\psi}}

\def\Phihat{{\widehat{\Phi}}}
\def\Phat{{\widehat{\Phi}}}

\def\K{{{\cal K}_\infty}}

\def\phat{\widehat{\ph}}

\def\Pepsi{\hbox{P\'e}^*}

\def\Wbar{\overline{W}}

\def\G{{\cal G}}





\topmatter
\title
INTEGRAL MOMENTS OF AUTOMORPHIC L--FUNCTIONS
\endtitle
\author
Adrian Diaconu\\ Paul Garrett
\endauthor
\address
Adrian Diaconu, School of Mathematics, University of Minnesota, Minneapolis, MN 55455
\endaddress
\email cad\@math.umn.edu
\endemail

\address
Paul Garrett, School of Mathematics, University of Minnesota, Minneapolis, MN 55455
\endaddress
\email garrett\@math.umn.edu
\endemail

\abstract This paper exposes the underlying mechanism for obtaining
second integral moments of $GL_2$ automorphic $L$--functions over an
arbitrary number field. Here, moments for $GL_2$ are presented in a form enabling application
of the structure of adele groups and their representation theory. To
the best of our knowledge, this is the first formulation of integral
moments in adele-group-theoretic terms, distinguishing global and
local issues, and allowing uniform application to number fields.
When specialized to  
the field of rational numbers $\Bbb{Q}$, we recover the classical results. 
\endabstract
\endtopmatter

\def\hk{\hfill\break}

\noindent
1. Introduction \hk
2. Poincar\'e series \hk
3. Unwinding to Euler product \hk
4. Spectral decomposition of Poincar\'e series \hk
5. Asymptotic formula \hk
Appendix 1: Convergence of Poincar\'e series \hk
Appendix 2: Mellin transform of Eisenstein Whittaker functions

\vskip30pt

\noindent{\bf \S 1.  Introduction}
\vskip 10pt

\noindent For ninety years, the study of mean values of families of automorphic
$L$--functions  
has played a central role in analytic number theory, for
applications to classical problems. In the absence of the Riemann
Hypothesis, or the Grand Riemann Hypothesis, rather, when referring
to general $L$--functions, suitable mean value results often serve as
a substitute. In particular, obtaining asymptotics or
sharp bounds for integral moments of automorphic
$L$--functions is of considerable interest. The study of integral
moments was initiated in 1918 by Hardy and Littlewood (see
\cite{Ha-Li}) who obtained the asymptotic formula for second moment
of the Riemann zeta-function  
$$\int_0^T \left|\zeta\left({\scriptstyle \frac12 }+ it\right)\right|^2 \, dt \; \sim \; T\log T \tag 1.1$$ 
About 8 years later, Ingham in \cite{I} obtained the fourth moment
 $$\int_0^T \left|\zeta\left({\scriptstyle \frac12 }+ it\right)\right|^4 \, dt \; \sim \; \frac{1}{2\pi^2}\cdot T(\log T)^4 \tag 1.2$$ 
 Since then, many papers by various authors have been devoted to this
 subject. For instance 
 see  \cite{At}, \cite{H-B}, \cite{G1}, \cite{M1}, \cite{J1}. Most
 existing results concern integral  
 moments of automorphic $L$--functions for $GL_1(\Q)$ and $GL_2(\Q)$. No analogue of (1.1) or (1.2) is known over an arbitrary number field. The only previously known results, for fields other than $\Bbb{Q}$, 
 are in \cite{M4}, \cite{S1}, \cite{BM1}, \cite{BM2} and \cite{DG2}, all over quadratic fields.

Here we expose the underlying mechanism to obtain second integral
moments of $GL_2$ automorphic $L$--functions over an arbitrary number
field. Integral moments for $GL_2$ are presented in a form amenable to
application of the representation theory of adele groups. To the best
of our knowledge, this is the first formulation of integral moments on
adele groups, distinguishing global and local questions, and allowing
uniform application to number fields.
More precisely, for $f$ an automorphic form on $GL_2$ and $\chi$ an
idele class character of the number field, let $L(s, f\otimes \chi)$
denote the twisted $L$--function attached to $f$. We obtain
asymptotics for averages 
$$\sum_{\chi}\;\,\int\limits_{-\infty}^{\infty}\;
\left|L\big({\scriptstyle \frac{1}{2}} + it,\, f\otimes
\chi\big)\right|^2 M_{\chi}(t)\, dt \tag 1.3$$ for suitable smooth
weights $M_{\chi}(t)$. The sum in (1.3) is over a certain set of idele
class characters which is infinite, in general. For general number
fields, it seems that (1.3) is the correct structure of the second
integral moment of $GL_2$ automorphic $L$--functions. This was first
pointed 
out  by  
Sarnak in \cite{S1}, where an average of the above type was studied
over the Gaussian field   
$\Bbb{Q}(i);$ see also \cite{DG2}. From the analysis of Section 2, it will become apparent that this 
comes from the Fourier transform on the idele class group of the
field.         

Meanwhile, in joint work \cite{DGG} with Goldfeld, the present authors
have found an extension to treat integral moments for $GL_r$ over
number fields. We exhibit specific Poincar\'e series $\Pe$ giving
identities of the form
$$ 
\hbox{moment expansion} \,= \int_{Z_\A GL_r(k)\\ GL_r(\A)} \Pe\cdot 
|f|^2 
\,= \,
\hbox{spectral expansion} 
$$ 
for cuspforms $f$ on $GL_r$.
The moment expansion on the left-hand side is of the form                               
$$ 
\sum_{F} {1\over 2\pi i} \int\limits_{\Re(s)={1\over 2}} |L(s,f\otimes F)|^2\;M_F(s)\;ds \;+\; \ldots
$$ 
summed over $F$ in an orthonormal basis for cuspforms on $GL_{r-1}$, as
well as corresponding continuous-spectrum terms. The specific choice
gives a kernel with a surprisingly simple {\it spectral expansion},
with only three parts: a leading term, a sum induced from cuspforms on
$GL_2$, and a continuous part again induced from $GL_2$. In
particular, no cuspforms on $GL_\ell$ with $2<\ell\le r$
contribute. 
Since the discussion for $GL_r$ with $r>2$ depends essentially on the
details of the $GL_2$ results, the $GL_2$ case merits special
attention. We give complete details for $GL_2$ here. For $GL_2$ over
$\Bbb{Q}$ and square-free level, the average of moments has a single
term, recovering the classical integral moment   
$$\int\limits_{-\infty}^{\infty}\, \left| L\left({\scriptstyle\frac{1}{2}} + it,
f\right)\right|^{2}M(t)\, dt$$  

As a non-trivial example, consider the case of a cuspform $f$
on $GL_3$ over $\Q$. We construct a weight function 
$\Gamma(s, w, f_\infty, F_\infty)$ depending upon complex parameters 
$s$ and $w,$ and upon the {\it archimedean} data for both $f$ and 
cuspforms $F$ on $GL_2,$ such that $\Gamma(s, w, f_\infty, F_\infty)$ 
has explicit asymptotic behavior similar to those in Section 5 below,
and such that the {\it moment expansion} above becomes
$$ 
\int_{Z_\A GL_3(\Q) \\ GL_3(\A)} \Pe(g) \, |f(g)|^2\,dg\;\;
= 
\sum_{F\;\text {on}\;GL_2} 
{1\over 2\pi i} 
\int\limits _{\Re(s)=\frac{1}{2}} 
|L(s,f\otimes F)|^2\cdot\Gamma(s, w, f_\infty, F_\infty)\,ds 
$$ 
$$ 
+ \,\,
{1\over 4\pi i}{1\over 2\pi i} 
\sum_{k\in {\Bbb Z}}
\;\,
\int\limits_{\Re(s_1)=\frac{1}{2}} \;
\int\limits_{\Re(s_2)=\frac{1}{2}} \;
\!\!\!\!\!|L(s_1, f\otimes E^{(k)}_{1-s_2})|^2  
\cdot 
\Gamma(s_1, w, f_\infty, E^{(k)}_{1-s_2,\infty})\;ds_2\,ds_1 
$$
where
$$ 
L(s_1, f\otimes E_{1-s_2}^{(k)}) 
= 
{L(s_1 - s_2 + {\scriptstyle\frac{1}{2}}, f) \cdot L(s_1 + s_2 - {\scriptstyle\frac{1}{2}} ,  f) 
\over 
\zeta(2-2s_2) 
} 
$$ 
In the above 
expression, $F$ runs over an orthonormal basis for all level-one 
cuspforms on $GL_2,$ with {\it no} restriction on the right 
$K_\infty$--type. Similarly, the Eisenstein series $E_s^{(k)}$ run over all
level-one Eisenstein series for $GL_2(\Q)$ with no restriction on
$K_\infty$--type, denoted here by $k.$


The course of the argument makes several points clear. First, the sum
of moments of twists of $L$--functions has a natural integral
representation. Second, the kernel arises from a collection of {\it
local} data, wound up into an automorphic form, and the computation
proceeds by unwinding. Third, the local data at {\it finite} primes is of a
mundane sort, already familiar from other constructions. Fourth, the
only subtlety resides in choices of archimedean data. Once this is
understood, it is clear that Good's  
original idea in \cite{G2}, seemingly limited to $GL_2(\Bbb{Q})$, exhibits a 
good choice of local data for {\it real} primes. See also \cite{DG1}. 
Similarly, while \cite{DG2} explicitly addresses only $GL_2(\Bbb{Z}[i])$, the 
discussion there exhibits a good choice of local data for {\it 
complex} primes. That is, these two examples suffice to illustrate the 
non-obvious choices of local data for all archimedean places. 

The structure of the paper is as follows. In Section 2, a family of
Poincar\'e series is defined in terms of local data, abstracting
classical examples in a form applicable to $GL_r$ over a number field.
In Section 3, the integral of the Poincar\'e series against $|f|^2$
for a cuspform $f$ on $GL_2$ is unwound and expanded, yielding a sum
of weighted moment integrals of $L$-functions $L(s,f\otimes \chi)$ of
twists of $f$ by Gr\"o\ss encharakteren $\chi$.  In Section 4, we find the
spectral decomposition of the Poincar\'e series: the leading term is
an Eisenstein series, and there are cuspidal and continuous-spectrum
parts with explicit coefficients.  In section 5, we derive an
asymptotic formula for integral moments, and observe that the length
of the averages involved is suitable for subsequent applications to
convexity breaking in the $t$--aspect.  The first appendix discusses
convergence of the Poincar\'e series in some detail, proving pointwise
convergence from two viewpoints, also proving $L^2$ convergence. The
second appendix computes integral transforms necessary to understand
the details in the spectral expansion.

For applications, one needs to combine refined choices of archimedean
data with extensions of the estimates in \cite{Ho-Lo} and \cite{S2}
(or \cite{BR}) to number fields. However, for now, we content
ourselves with a formulation that lays the groundwork for applications
and extensions.  In subsequent papers we will address convexity
breaking in the $t$--aspect, and extend this approach to $GL_r$.



\vfill\break

\vskip20pt\noindent {\bf \S 2. Poincar\'e series}
\vskip 10pt 

Before introducing our Poincar\'e series $\Pe(g)$ on $GL_{r}$ ($r \ge 2$) mentioned 
in the introduction, we find it convenient 
to first fix some notation in this context. 
Let $k$ be a number field, $G=GL_r$ over $k,$ and define the standard 
subgroups:
$$ 
P = P^{r-1,1} 
= \left\{\pmatrix \hbox{$(r-1)$-by-$(r-1)$} & * \cr 0 & \hbox{$1$-by-$1$}\endpmatrix \right\} 
$$ 
the standard maximal proper parabolic subgroup,   
$$ 
U = \left\{\pmatrix I_{r-1} & * \cr 0 & 1 \endpmatrix\right\} 
\qquad  
H = \left\{\pmatrix \hbox{$(r-1)$-by-$(r-1)$} & 0 \cr 0 & 1 \endpmatrix \right \} 
\qquad 
Z = \hbox{center of $G$}
$$
Let $K_{\nu}$ denote the standard maximal 
compact in the $k_{\nu}$--valued points $G_{\nu}$ of $G.$

The {\bf Poincar\'e series} $\Pe(g)$ is of the form
$$ 
\Pe(g)\;\;\; = \sum_{\gamma\in Z_k H_k\\G_k} \ph(\gamma g) \qquad \qquad (g\in G_{\A}) \tag 2.1
$$
for suitable functions $\ph$ on $G_{\A}$ described as follows. For $v \in \Bbb{C},$ let 
$$ 
\ph \,=\, \bigotimes_{\nu} \, \ph_{\nu} \tag 2.2
$$ 
where for $\nu$ {\it finite} 
$$ 
\ph_{\nu}(g)\, 
=\, \cases \left|(\text{det} \, A)/d^{r - 1}\right|_{\nu}^{v}  & \text{for $g=mk$ with $m = \pmatrix A & 0\cr 0 & d \endpmatrix \in Z_{\nu} H_{\nu}$ and $k\in K_{\nu}$}\cr
0  & \text{otherwise} \endcases  \tag 2.3  
$$ 
and for $\nu$ {\it archimedean} require right $K_{\nu}$--invariance and left 
equivariance 
$$ 
\ph_{\nu}(mg) 
= \left|{\text{det}\,  A\over d^{r-1}}\right|_{\nu}^{v} \cdot \ph_{\nu}(g) \qquad 
\left(\text{for $g\in G_{\nu}$ and $m = \pmatrix A & 0\cr 0 & d \endpmatrix \in Z_{\nu} H_{\nu}$} \right) \tag 2.4
$$
Thus, for $\nu |\infty,$ the further data determining $\ph_{\nu}$ consists of 
its values on $U_{\nu}.$ The simplest useful choice is 
$$ 
\ph_{\nu}\pmatrix I_{r-1} & x \cr 0 & 1 \endpmatrix = \left(1+ |x_1|^2 + \cdots + 
|x_{r-1}|^2 \right)^{- d_{\nu} (r - 1)w_{\nu}/2} \;\;\;\; 
\left(\text{$x=\pmatrix x_1\cr \vdots\cr x_{r-1} \endpmatrix$ and $w_{\nu} \in \C$}\right)  \tag 2.5
$$ 
with $d_{\nu} = [k_\nu:\R]$. Here 
the norm $|x_1|^2+\cdots+|x_{r-1}|^2$ is 
invariant under $K_{\nu},$ that is, $|\cdot|$ is the usual absolute value on 
$\Bbb{R}$ or $\Bbb{C}.$ Note that by the product formula $\ph$ is left 
$Z_{\A} H_{k}$--invariant.

We have the following

\vskip 12pt 
\proclaim{Proposition 2.6}{\rm (Apocryphal)} With the specific choice (2.5)
of $\ph_\infty = \otimes_{\v|\infty}\,\ph_\v$, the series (2.1) defining $\Pe(g)$ 
converges absolutely and locally uniformly for $\Re(\sv)>1$ and $\Re(\sw_\v)>1$ for all
$\v|\infty$.  
\endproclaim

\vskip10pt
\noindent {\bf Proof:} In fact, the argument applies to a much broader
class of archimedean data. For a complete argument when $r = 2$, 
and $\sw_{\nu} =\sw$ for all $\nu| \infty$, see Appendix 1. \qquad
\qed

We can give a broader and more robust, though somewhat weaker, result, 
as follows. Again, for simplicity, we shall assume $r = 2$. Given 
$\ph_\infty$, for $x$ in 
$k_\infty=\prod_{\v|\infty}k_\v$, let
$$
\Phi_\infty(x) \,=\, \ph_\infty\bmatrix{1&x\cr 0&1}
$$
For $0<\ell\in\Z$, let $\Omega_\ell$ be the collection of $\ph_\infty$ 
such that the associated $\Phi_\infty$ is absolutely integrable, and such
that the Fourier transform $\Phihat_\infty$ along $k_\infty$ satisfies
the bound 
$$
\Phihat_\infty(x) \;\ll\; \underset{\nu|\infty}\to{\textstyle\prod}\, (1 + |x|_{\nu}^2)^{-\ell}
$$ 
For example, for $\ph_\infty$ to be in $\Omega_\ell$ it suffices that
$\Phi_\infty$ is $\ell$ times continuously differentiable, with each 
derivative absolutely
integrable. For $\Re(\sw_{\nu}) > 1$, 
$\nu|\infty$, the simple explicit
choice of $\ph_\infty$ above lies in $\Omega_\ell$ for {\it every}
$\ell>0$.

\vskip 12pt
\proclaim{Theorem 2.7}{\rm (Apocryphal)} 
Suppose $r = 2$, $\Re(v),\, \ell$ sufficiently large, and $\ph_\infty \in \Omega_\ell$. 
The series defining $\Pe(g)$ converges absolutely and
locally uniformly in both $g$ and $v$. Furthermore, up to an Eisenstein series, 
the Poincar\'e series is square integrable on $Z_\A G_k\\G_\A$.  
\endproclaim

\vskip10pt
\noindent{\bf Proof:} See Appendix 1. \qquad
\qed

The precise Eisenstein series to be subtracted from the Poincar\'e
series to make the latter square-integrable will be discussed in
Section 4 (see formula 4.6). For our special choice (2.5) of
archimedean data, both these convergence results apply with $\Re(\sw_\nu)
> 1$ for $\nu|\infty$ and $\Re(\sv)$ large.

For convenience, a monomial vector $\ph$ as in (2.2) described by 
(2.3) and (2.4) will be called {\it admissible}, if $\ph_\infty \in \Omega_\ell$, 
with both $\Re(v)$ and $\ell$ sufficiently large. 

\vskip20pt 
\noindent{\bf \S 3. Unwinding to Euler product}
\vskip 10pt 

From now on, we shall assume $r = 2.$ Recall the notation made in the previous 
section, which in the present case reduces to: $G = GL_{2}$ over the number field $k$ 
together with the standard subgroups    
$$ P = \left\{\pmatrix*&*\cr 0&* \endpmatrix \right\} 
\;\;\;\;\; 
N = U = \left\{\pmatrix1&*\cr 0&1\endpmatrix \right\} 
\;\;\;\;\; 
M = ZH = \left\{\pmatrix*&0\cr 0&* \endpmatrix \right\} 
$$ 
Also, for any place $\nu$ of $k$, let $K_{\nu}$ be the standard maximal 
compact subgroup. That is, for finite $\nu,$ we take 
$K_{\nu} = GL_2(\o_{\nu}),$ at real places $K_{\nu}=O(2),$ and at complex 
places $K_{\nu} = U(2).$

With the Poincar\'e series defined by (2.1), our main goal is to 
unwind a corresponding global integral to express it as an inverse 
Mellin transform of an Euler product. For convenience, recall that   
$$\Pe(g) \;\;= \sum_{\gam\in M_k \backslash G_k} \ph(\gam g) 
\qquad \qquad (g\in G_\A) 
\tag 3.1$$ 
where the {\it monomial} vector 
$$\ph = \bigotimes_{\nu} \ph_{\nu}
$$ 
is defined by 
$$\ph_\nu(g)\, 
=\, \cases \chi_{0,\nu}(m)  & \text{for $g=mk$, $m\in M_\nu$ and $k\in K_{\nu}$}\cr
0  & \text{for $g\not\in M_\nu\cdot K_{\nu}$} \endcases   
\qquad \;\;\text{(for $\nu$ finite)}
\tag 3.2$$ and for $\nu$ infinite, we do not entirely specify 
$\ph_{\nu},$ only requiring the left equivariance 
$$\ph_\nu(m n k) = \chi_{0,\nu}(m)\cdot \ph_\nu(n) 
\;\;\;\;\;\;\; \text{(for $\nu$ infinite, $m\in M_\nu$, $n\in N_\nu$ and $k\in K_{\nu}$)} \tag 3.3$$ 
Here, $\chi_{0,\nu}$ is the character of $M_{\nu}$ given by 
$$\chi_{0,\nu}(m) = \left| \frac{a}{d} \right|_{\nu}^{v} \qquad \qquad \left(m =  \pmatrix a & 0 \cr 0 & d \endpmatrix \in M_{\nu},\, v\in \Bbb C \right) \tag 3.4$$ Then, 
$\chi_0=\bigotimes_\nu \chi_{0,\nu}$ is  
$M_k$--invariant, and $\ph$ has trivial central character and is left
$M_\A$--equivariant by $\chi_0.$ Also, note that for $\nu$
infinite, our assumptions imply that 
$$x \to \ph_\nu\pmatrix1& x \cr 0&1\endpmatrix$$ 
is a function of $|x|$ only.

Let $f_1$ and $f_2$ be cuspforms on $G_\A$. Eventually we will take
$f_1=f_2$, but for now merely require the following. At all
$\nu$, require (without loss of generality) that $f_1$ and
$f_2$ have the same right $K_\nu$--type, that this $K_\nu$--type is
{\it irreducible}, and that $f_1$ and $f_2$ correspond to the same
vector in the $K$--type (up to scalar multiples). Schur's lemma
assures that this makes sense, insofar as there are no non-scalar
automorphisms. Suppose that the representations of $G_\A$ 
generated by $f_1$ and $f_2$ are {\it irreducible}, with the same
central character. Last, require that each $f_i$ is a
special vector locally everywhere in the representation it generates, 
in the following sense. Let
$$f_i(g) \;\;= \sum_{\xi\in Z_k\\M_k} W_i(\xi g)\tag 3.5 
$$
be the Fourier expansion of $f_i$, and let 
$$
W_i=\bigotimes_{\nu\le \infty} W_{i,\nu}
$$
be the factorization of the Whittaker function $W_i$ into local
data. By \cite{JL}, we may require that for all $\nu<\infty$ the Hecke type
local integrals
$$
\int\limits_{a\in k_{\nu}^\times} 
W_{i,\, \nu} 
\pmatrix a & 0 \cr 0 & 1\endpmatrix
\,|a|_{\nu}^{s-\frac{1}{2}}\,da
$$
differ by at most an exponential function from the local $L$--factors
for the representation generated by $f_i$. Eventually we will take
$f_1=f_2$, compatible with these requirements.

The integral under consideration is (with notation suppressing 
details)
$$I(\chi_0) 
= \int_{Z_\A G_k \\ G_\A} \Pe(g)\,f_1(g)\,\bar {f}_2(g)\,dg\tag 3.6 
$$ 

For $\chi_0$ (and archimedean data) in the range of absolute 
convergence, the integral unwinds (via the definition of 
the Poincar\'e series) to 
$$ \int_{Z_\A M_k \\ G_\A} \ph(g)\,f_1(g)\,\bar{f}_2(g)\,dg
$$ 
Using the Fourier expansion 
$$ f_1(g) \;\,= \sum_{\xi\in Z_k\\M_k} W_1(\xi\,g) 
$$ 
this further unwinds to 
$$ \int_{Z_\A \\ G_\A} \ph(g)\,W_1(g)\,\bar{f}_2(g)\,dg \tag 3.7 
$$ 

Let $C$ be the idele class group $GL_1(\A)/GL_1(k)$, and $\Chat$ its 
dual. More explicitly, by Fujisaki's Lemma (see Weil \cite{W1}, page 32, Lemma 3.1.1), 
the idele class group $C$ is a product of a copy of $\R^+$ 
and a compact group $C_0$. By Pontryagin duality, $\Chat \isom \R \times \Chat_0$ with 
$\Chat_0$ discrete. It is well-known that, for any compact open subgroup $U_{\text{fin}}$ of 
the finite-prime part in $C_0$, the dual of $C_{0}/ U_{\text{fin}}$ is
finitely generated with rank $[k:\Bbb{Q}] - 1$. 
The general Mellin transform and inversion are  
$$\align f(x) \,&= \int_\Chat \int_C f(y) \chi(y)\,dy\, \chi^{-1}(x)\, d\chi \tag 3.8\\&
=  \sum_{\chi' \in \Chat_0}\;\; \frac{1}{2\pi i}\int\limits_{_{\Re(s) =\sigma}} \int_C f(y) \chi'(y)|y|^{s}\,dy\, {\chi'}^{-1}(x)|x|^{-s}\, ds
\endalign$$ 
for a suitable Haar measure on $C$. 

To formulate the main result of this section, we need one more piece of notation. For $\nu$ infinite and $s\in \Bbb{C}$, let  
$$\align \Cal{K}_{\nu}(s,\, \chi_{0, \nu},\, \chi_{\nu}) \, & =  \,\int_{Z_{\nu}\backslash M_{\nu}N_{\nu}} \int_{Z_{\nu}\backslash M_{\nu}}
\ph_{\nu}(m_{\nu}n_{\nu})W_{1, \nu}(m_{\nu}n_{\nu})\\
& \hskip 5pt\cdot \Wbar_{2, \nu}({m}_{\nu}'n_{\nu})\,\chi_{\nu}({m}_{\nu}')\,|m_{\nu}'|_{\nu}^{s - \frac{1}{2}}\,\chi_{\nu}(m_{\nu})^{-1}\, |m_{\nu}|_{\nu}^{\frac{1}{2} - s}\,d{m}_{\nu}'\,dn_{\nu} \, dm_{\nu}\tag 3.9\endalign$$ and set 
$$\Cal{K}_{\infty}(s,\, \chi_{0},\, \chi) \,
=\, \prod_{\nu |\infty} \, \Cal{K}_{\nu}(s,\, \chi_{0, \nu},\,
\chi_{\nu})\tag 3.10$$ Here $\chi_0=\bigotimes_\nu \chi_{0,\nu}$ is
the character defining the monomial vector $\ph$, and
$\chi=\bigotimes_\nu \chi_{\nu}\in \Chat_0$. 
When the monomial vector $\ph$ is {\it admissible}, the integral (3.9)
defining $K_\nu$ converges absolutely for $\Re(s)$ sufficiently large.
We are especially interested in the choice
$$\ph_{\nu}(n)\, 
=\, \cases \left(1 + x^2\right)^{-\frac{w}{2}}  & \text{for $\nu |\infty$ real, and $n = \pmatrix 1 & x \cr 0& 1 \endpmatrix  \in N_{\nu}$}\\
\left(1 + |x|^{2} \right)^{- w}  & \text{for $\nu |\infty$ complex, and $n = \pmatrix 1 & x \cr 0& 1 \endpmatrix \in N_{\nu}$} \endcases \;\;\;\;\; (v,\, w\in \Bbb C)
\tag 3.11$$
The monomial vector $\ph$ generated by this choice is admissible for
$\Re(w)>1$ and $\Re(v)$ sufficiently large. This choice will be used in Section 5 to
derive an asymptotic formula for the $GL_2$ integral moment over the
number field $k$. The main result of this section is  
 
\vskip 12pt
\proclaim{Theorem 3.12} 
For $\ph$ an admissible monomial vector as above, for suitable
$\sigma>0$, 
$$
I(\chi_0)  \; =  \sum_{\chi \in \Chat_0}\;\; \frac{1}{2\pi i}\int\limits_{_{\Re(s) =\sigma}} L(\chi_{0} \cdot \chi^{-1}| \cdot |^{1 - s},\,f_1) 
\cdot 
L(\chi |\cdot|^{s},\,\bar{f}_{2})\, 
\Cal{K}_{\infty}(s,\, \chi_{0},\, \chi) \, ds$$ 
Let $S$ be a finite set of places including archimedean places,
all absolutely ramified primes, and all finite bad places for
$f_1$ and $f_2$. Then the sum is over a set $\Chat_{0,S}$ of
characters unramified outside $S$, with bounded ramification at
finite places, depending only upon $f_1$ and $f_2$.
\endproclaim

\vskip 10pt 
\noindent{\bf Proof:} Applying (3.8) to 
$\bar{f}_2$ via the identification $$\left\{\pmatrix a'&0\cr
0&1 \endpmatrix :  a' \in C\right\} \isom C$$  
and using the Fourier expansion 
$$f_2(g) \;\,= \sum_{\xi\in Z_k\\M_k} W_2(\xi\,g) 
$$ 
the integral (3.7) is 

\vbox{ 
$$ \int_{Z_\A \\ G_\A} \ph(g)\,W_1(g)\, 
\left(\int_\Chat \int_{C} 
\bar{f}_2(m'g)\,\chi(m')\,dm'\,d\chi\right)dg 
$$ 
$$ 
= \int_\Chat\left( 
\int_{Z_\A \\ G_\A} \ph(g)\,W_1(g)\, 
\int_{C}\; 
\sum_{\xi\in Z_k\\M_k} \Wbar_2(\xi m'g)\,\chi(m')\,dm'\,dg 
\right)d\chi 
$$ 
$$ 
= \int_\Chat\left( 
\int_{Z_\A \\ G_\A} \ph(g)\,W_1(g)\, 
\int_{\J} \;
\Wbar_2(m'g)\,\chi(m')\,dm'\,dg 
\right)d\chi
$$ 
} 
\noindent where $\J$ is the ideles. The interchange of order of
integration is justified by the absolute convergence of
the outer two integrals. (The innermost integral {\it
cannot} be moved outside.) This follows from the rapid decay of
cuspforms along the split torus.

For fixed $f_1$ and $f_2$, the finite-prime ramification of the 
characters $\chi\in \Chat$ is bounded, so there are only finitely many 
bad finite primes for all the $\chi$ which appear. In particular, all the 
characters $\chi$ which appear are unramified outside $S$ and with bounded 
ramification, depending only on $f_1$ and $f_2$, at finite places in $S$. Thus, 
for $\nu \in S$ finite, there exists a compact open subgroup $U_{\nu}$ of $\o_{\nu}^{\times}$ 
such that the kernel of the $\nu^{\text{th}}$ component $\chi_{\nu}$ of $\chi$ contains $U_{\nu}$ 
for all characters $\chi$ which appear.

Since
$f_1$ and $f_2$ generate irreducibles
locally everywhere, the Whittaker functions $W_i$ 
factor 
$$
W_i(\{g_\nu:\nu\le \infty\}) = \Pi_\nu W_{i,\nu}(g_\nu) 
$$ 
Therefore, the inner integral over $Z_\A \\G_\A$ and $\J$ 
factors over primes, and 
$$I(\chi_0) 
= \int_\Chat 
\Pi_\nu 
\left( 
\int_{Z_\nu \\ G_\nu} 
\int_{k_{\nu}^{\times} } 
\ph_\nu(g_\nu)\,W_{1,\nu}(g_\nu)\, 
\Wbar_{2,\nu}(m_\nu'g_\nu)\,\chi_\nu(m'_\nu)\,dm_\nu'\,dg_\nu 
\right)d\chi 
$$ 

Let $\omega_{\nu}$ be the $\nu^{\text{th}}$ component of the central character 
$\omega$ of $f_2$. Define a character of $M_\nu$ by 
$$\pmatrix a &0\cr 0& d \endpmatrix \to \omega_{\nu}\pmatrix d
&0\cr 0& d \endpmatrix \chi_{\nu}\pmatrix a/d &0\cr 0& 1
\endpmatrix$$ Still denote this character by $\chi_{\nu},$ without
danger of confusion. In this notation, the last 
expression of $I(\chi_0)$ is
$$I(\chi_0) 
= \int_\Chat 
\Pi_\nu 
\left( 
\int_{Z_\nu \\ G_\nu} 
\int_{Z_{\nu} \\ M_{\nu} } 
\ph_\nu(g_\nu)\,W_{1,\nu}(g_\nu)\, 
\Wbar_{2,\nu}(m_\nu'g_\nu)\,\chi_\nu(m'_\nu)\,dm_\nu'\,dg_\nu 
\right)d\chi 
$$

Suppressing the index $\nu$, the $\nu^{\text{th}}$ local integral is 
$$ 
\int_{Z\\G} 
\int_{Z\\M} 
\ph(g)\,W_1(g)\, 
\Wbar_2(m'g)\,\chi(m')\,dm'\,dg 
$$ 

Take $\nu$ {\it finite} such that both $f_1$ and $f_2$ are right
$K_{\nu}$--invariant. Use a $\nu$--adic Iwasawa decomposition
$g=mnk$ with $m\in M,$ $n\in N,$ and $k\in K.$ The Haar
measure is $d(mnk)=dm\,dn\,dk$ with Haar measures on the factors.
The integral becomes
$$
\int_{Z\\MN} 
 \int_{Z\\M}
\ph(mn)\,W_1(mn)\,
\Wbar_2(m'mn)\,\chi(m')\,dm'\,dn\,dm
$$
To symmetrize the integral, replace $m'$ by $m'm^{-1}$ to obtain
$$
\int_{Z\\MN} \int_{Z\\M}
\ph(mn)\,W_1(mn)\,
\Wbar_2(m'n)\,\chi(m')\,\chi(m)^{-1}\,dm'\,dn\,dm
$$
The Whittaker functions $W_i$ have left $N$--equivariance
$$W_i(ng) = \psi(n)\,W_i(g)
\;\;\;\;\hbox{(fixed non-trivial $\psi$)}
$$
so
$$ W_1(mn) =  W_1(mnm^{-1}\,m) 
= \psi(mnm^{-1})\,W_1(m)
$$
and similarly for $W_2$. Thus, letting
$$ X(m,m') = 
\int_N
\ph(n)\,
\psi(mnm^{-1})\,
\overline{\psi}(m'n{m'}^{-1})\,
dn
$$
the local integral is
$$
\int_{Z\\M} 
 \int_{Z\\M}
\chi_0(m)\,
W_1(m)\,
\Wbar_2(m')\,\chi(m')\,\chi^{-1}(m)\,
X(m,m')
\,dm'\,dm
$$
We claim that for $m$ and $m'$ in the supports of the Whittaker
functions, the inner integral $X(m,m')$ is constant, independent of
$m,\, m',$ and it is $1$ for almost all finite primes. First, $\ph(mn)$
is $0,$ unless $n\in M\cdot K \cap N,$ that is, unless $n\in
N \cap K.$ On the other hand, 
$$
\psi(mnm^{-1})\cdot W_1(mk) 
= \psi(mnm^{-1})\cdot W_1(m) 
= W_1(mn)
= W_1(m)
\;\;\;\;\;\;\hbox{(for $n\in N\cap K$)}
$$
Thus, for $W_1(m)\not=0$, necessarily $\psi(mnm^{-1})=1$. A similar
discussion applies to $W_{2}$. 
So, up to normalization, the inner integral is $1$ for
$m,\, m'$ in the supports of $W_{1}$ and $W_{2}$. Then 

\vbox{
$$
\int_{Z\\M} 
 \int_{Z\\M}
\chi_0(m)\,
W_1(m)\,
\Wbar_2(m')\,\chi(m')\,\chi^{-1}(m)
\,dm\,dm'
$$
$$
=
\int_{Z\\M}
(\chi_0\cdot \chi^{-1})(m)
\, W_1(m)\,
dm
\cdot
\int_{Z\\M}
\chi(m')\,
\Wbar_2(m')\,
dm'
$$
$$ 
= L_\nu(\chi_{0, \nu} \cdot \chi_{\nu}^{-1}|\cdot|_{\nu}^{1/2},\,f_1) 
\cdot 
L_\nu(\chi_{\nu}|\cdot|_{\nu}^{1/2},\,\bar{f}_2)
$$ 
} 
\noindent i.e., the product of local factors of the standard
$L$--functions in the theorem (up to exponential functions at
finitely many finite primes) by our assumptions on $f_1$ and $f_2$.

For non-trivial right $K$--type $\sigma$, the argument is similar but a
little more complicated. The key point is that the inner integral over
$N$ (as above) should not depend on $mk$ and $m'k$, for $mk$ and $m'k$
in the support of the Whittaker functions. Changing conventions
for a moment, look at $V_\sigma$--valued Whittaker functions, and
consider any $W$ in the $\nu\th$ Whittaker space for $f_i$ having
right $K$--isotype $\sigma$. Thus,
$$W(gk) = \sigma(k)\cdot W(g)
\;\;\;\;\;\text{(for $g\in G$ and $k\in K$)}
$$
For $\ph(mn)\not=0$, again $n\in N\cap K$. Then
$$ \sigma(k)\cdot \psi(mnm^{-1})\cdot W(m) 
= W(mnk) 
= \sigma(k)\cdot W(mn) 
= \sigma(k)\cdot \sigma(n)\cdot W(m) 
$$
where in the last expression $n$ comes out on the right by the right
$\sigma$--equivariance of $W$. For $m$ in the support of $W$, 
$\sigma(n)$ acts by the scalar $\psi(mnm^{-1})$ on $W(mk)$, for all
$k\in K$. Thus, 
$\sigma(n)$ is scalar on that copy of $V_\sigma$. At the same time,
this scalar is $\sigma(n)$, so is independent of $m$ if
$W(m)\not=0$. Thus, except for a common integral over $K$, the
local integral falls into two pieces, each yielding the local
factor of the $L$--function. The common integral over $K$ is
a constant (from Schur orthogonality), non-zero since the two vectors
are collinear in the $K$--type. \qquad
\qed

\vskip10pt

At this point the archimedean local factors of the Euler 
product are not specified. The option to vary the choices is 
{\it essential} for applications.

\vskip20pt\noindent {\bf \S 4. Spectral decomposition of Poincar\'e series}
\vskip 10pt

The objective now is to spectrally decompose the Poincar\'e series
defined in (3.1). Throughout this section, we assume that $\ph$ is
{\it admissible}, in the sense given at the end of Section 2.  As we
shall see, in general $\Pe(g)$ is not square-integrable. However,
choosing the archimedean part of the monomial vector $\ph$ to have
enough decay, and after an obvious Eisenstein series is subtracted,
the Poincar\'e series is not only in $L^2$ but also has sufficient
decay so that its integrals against Eisenstein series converge
absolutely, by explicit computation. In particular, if the archimedean
data is specialized to (3.11), the Poincar\'e series $\Pe(g)$ has
meromorphic continuation in the variables $v$ and $w$. This is
achieved via spectral decomposition and meromorphic continuation of
the spectral fragments. See \cite{DG1}, \cite{DG2} when $k = \Bbb{Q}$,
$\Bbb{Q}(i)$. 

Let $k$ be a number field, $G=GL_2$ over $k,$ and $\omega$ a unitary character 
of $Z_k \\ Z_\A.$ Recall the decomposition 
$$L^2(Z_\A G_k\\G_\A,\, \omega)\;=\; L_{\text{cusp}}^2(Z_\A
G_k\\G_\A,\, \omega) \; \oplus \; L_{\text{cusp}}^2(Z_\A
G_k\\G_\A,\, \omega)^{\perp}$$ The orthogonal complement  
$$\align L_{\text{cusp}}^2(Z_\A G_k \backslash G_\A,\, \omega)^{\perp}\; &\isom\; \{1-\text{dimensional representations}\}\\
&\hskip15pt \oplus\, \int_{(GL_{1}(k)\backslash GL_{1}(\A))\,\widehat{\;}}^{\oplus}\;\; \underset{\nu}\to \bigotimes\; \text{Ind}_{P_{\nu}}^{G_{\nu}} (\chi_\nu \, \delta_{\nu}^{1/2})\, d\chi \endalign$$ 
where $\delta$ is the modular function on 
$P_\A,$ and the isomorphism is via Eisenstein series. Using this, we shall explicitly decompose 
our Poincar\'e series as 
$$
\Pe \,=\, \text{Eisenstein series \, $+$ \, discrete part \,$+$\,\,\,   continuous part} \qquad 
(\text{with $\omega = 1$})
$$

The projection to cuspforms is straightforward componentwise. We have

\vskip 12pt
\proclaim{Proposition 4.1} Let $f$ be a cuspform on $G_{\A}$ generating a spherical representation locally everywhere, and suppose $f$ corresponds to a spherical vector everywhere locally. In the region of absolute convergence of the Poincar\'e series $\Pe(g)$, the integral 
$$\int_{Z_\A G_k\\ G_\A} \bar{f}(g)\,\Pe(g)\,dg 
$$ 
is an Euler product. At finite $\nu$, the corresponding local 
factors are, up to a constant depending on the set of absolutely ramified primes
in $k,$ $L_\nu(\chi_{0, \nu}\,|\cdot|_{\nu}^{1/2},\, \bar{f}\,).$
\endproclaim

\vskip 10pt 
\noindent{\bf Proof:} The computation uses the same facts as the 
Euler factorization in the previous section. Using the Fourier 
expansion 
$$ f(g) \;\,= \sum_{\xi\in Z_k\\M_k} W(\xi g) 
$$ 
unwind 
$$ \int_{Z_\A G_k\\ G_\A} \bar{f}(g)\,\Pe(g)\,dg 
\;=\; \int_{Z_\A M_k\\ G_\A} \,\sum_{\xi}\; \overline{W}(\xi g)\,\ph(g)\,dg 
\;= \; \int_{Z_\A \\ G_\A} \overline{W}(g)\,\ph(g)\,dg 
$$ 
$$ 
=\; \prod_\nu \left(\int_{Z_\nu \\ G_\nu} 
\overline{W}_\nu(g_\nu)\,\ph_\nu(g_\nu)\,dg_\nu\right) 
$$ 
where the local Whittaker functions at finite places are normalized as in \cite{JL} 
to give the {\it correct} local $L$--factors.  

At finite $\nu$, suppressing the subscript $\nu$, the integrand 
in the $\nu^{\text{th}}$ local integral is right $K_{\nu}$--invariant, so we can 
integrate over $MN$ with left Haar measure. The 
$\nu^{\text{th}}$ Euler factor is 
$$ \int_{Z\backslash M}\int_N \overline{W}(mn) \, \ph(mn) \,\,dn\,dm 
= 
\int_{Z\backslash M}\int_N \overline{\psi}(mnm^{-1})\, \overline{W}(m) \, \chi_0(m)\,\ph(n) \,dn\,dm 
$$ for all finite primes $\nu$.  
The integral over $n$ is 
$$\int_N \overline{\psi}(mnm^{-1})\, \ph(n)\, dn 
$$ 
For $\ph(n)$ to be non-zero requires $n$ to lie in $M\cdot K$, which further  
requires, as before, that $n\in N\cap K$. Again, $W(m)=0$ unless 
$$m(N\cap K)m^{-1}\subset N\cap K 
$$ 
The character $\psi$ is trivial on $N\cap K$. Thus, the integral over 
$N$ is really the integral of $1$ over $N\cap K$. Thus, 
at finite primes $\nu$, the local factor is 
$$ 
\int_{Z\\M} \overline{W}(m)  \, \chi_0(m) \,dm 
\,= \,
L_\nu(\chi_{0, \nu}\,|\cdot|_{\nu}^{1/2},\,\bar{f}\,)  \qquad \qquad \qed
$$ 

Of course, the spectral decomposition of a right $K_{\A}$--invariant automorphic form can only involve 
everywhere locally spherical cuspforms. 

Assume that $\ph$ is given by (3.11). Taking $\Re(v) > 1$ and $\Re(w) > 1$ 
to ensure by Proposition 2.6 absolute convergence of $\Pe(g),$ the
local integral in Proposition 4.1 at infinite $\nu$ is
$$\int_{Z_\nu \\ G_\nu} 
\overline{W}_\nu(g_\nu)\,\ph_\nu(g_\nu)\,dg_\nu \,= \, \Cal{G}_{\nu}({\scriptstyle \frac{1}{2}} + i\bar{\mu}_{_{f, \nu}}; v,
w) 
$$ 
where, up to a constant,
$$\Cal{G}_{\nu}(s;
v, w) = \pi ^{- v}\,\frac{\Gamma\big(\frac{v + 1- s}{2}\big)\Gamma\big(\frac{v + w - s }{2}\big)\Gamma\big(\frac{v + s}{2}\big)\Gamma\big(\frac{v + w + s  -
1}{2}\big)}{\Gamma \big(\frac{w}{2}\big)\Gamma\big(v + \frac{w}{2}\big)} \tag 4.2$$
for $\nu | \infty$ real, and 
$$\Cal{G}_{\nu}(s; v,
w) = (2\pi)^{- 2v}\,\frac{\Gamma(v + 1 - s) 
\Gamma (v + w - s) \Gamma(v +  s) \Gamma (v + w +  s  - 1)}{\Gamma(w)\Gamma(2v + w)} \tag 4.3$$ 
for $\nu | \infty$ complex. In the above expression $i\mu_{_{f, \nu}}$ and $-i\mu_{_{f, \nu}}$ are 
the local parameters of $f$ at $\nu.$ These expressions as ratios of
products of gamma functions are obtained by standard computations (see
\cite{DG1} and \cite{DG2}), from the normalizations
$$ 
W_{\nu} \pmatrix a&\cr &1\endpmatrix  \,=\, 
\cases  |a|^{1/2} K_{i \mu_{_{f, \nu}}}(2 \pi |a|) 
& \text{if $\nu \isom \R$}\\      
|a|\, K_{2 i \mu_{_{f, \nu}}}(4 \pi |a|) & \text{if $\nu \isom \C$}
\endcases 
$$ 
and invocation of local multiplicity-one of Whittaker models. Then,
with respect to an orthonormal basis $\{F\}$ of everywhere  
locally spherical cuspforms, it is natural to consider the spectral sum  
$$
\sum_{F} \; \bar{\rho}_{_F}\,  \Cal{G}_{_{F_{_\infty}}}(v,
w) \, L(v + {\scriptstyle \frac{1}{2}},\overline{F})\cdot F  
$$ 
where
$$ 
\Cal{G}_{_{F_{_\infty}}}(v,
w) \,=\, \prod_{\nu | \infty} \, \Cal{G}_{\nu}({\scriptstyle \frac{1}{2}} + i\bar{\mu}_{_{F, \nu}}; v,
w) 
$$ 
with $\Cal{G}_{\nu}$ defined in (4.2) and (4.3). Here we absorbed all the 
{\it ambiguous} constants 
at infinite places into $\bar{\rho}_{_F}.$ Traditionally, the constant $\rho_{_F}$ is 
denoted by $\rho_{_F}(1)$ being considered the first Fourier coefficient of $F.$ As 
mentioned at the beginning of this section, and as we shall shortly see, 
the Poincar\'e series $\Pe(g)$ is up to an Eisenstein series a 
square-integrable function. It will then be clear that the above spectral sum 
represents the discrete part of $\Pe.$

By considering the 
usual integral representation against an Eisenstein series of the {\it
completed} $GL_2\times GL_2$ Rankin-Selberg 
$L$--function $\Lambda(s, F
\otimes \bar{F})$ \cite{J} (for the general case $GL_{m}\times GL_{n},$ 
see the review of the
literature in \cite{CPS2}), and then taking the residue at $s=1,$ one obtains  
$$\text{non-zero constant} \,
=\, |\rho_{_F}|^{2} \cdot L_{\infty}(1, F \otimes \bar{F})\cdot \underset{s = 1}\to {\text{Res}}\, L(s, F \otimes \bar{F})  \qquad (\text{for $||F|| = 1$})
$$ 
the constant on the left being independent of $F$. The local factors of $L(s, F \otimes \bar{F})$ on the right obtained from the integral representation may differ from those of the {\it correct} convolution $L$--function obtained from the local theory at only the absolutely ramified primes in $k.$ Comparing the gamma factors $L_{\infty}(1, F \otimes \bar{F})$ with $\Cal{G}_{_{F_{_\infty}}}(v,
w),$ we deduce that $\bar{\rho}_{_F}\,\Cal{G}_{_{F_{_\infty}}}(v,
w)$ has exponential decay in the local parameters of $F$. Combining
this with standard estimates and the Weyl's Law \cite{LV} (see also
\cite{Do} for an upper bound, which suffices for us), it follows that
the above spectral sum is absolutely convergent for $(v, w)\in
\Bbb{C}^2$, apart from the poles of $\Cal{G}_{_{F_{_\infty}}}(v,w).$

For the remaining decomposition, subtract (as in \cite{DG1},
\cite{DG2}) a finite linear combination of Eisenstein series from the
Poincar\'e series,  
leaving a function in $L^2$ with sufficient decay to be integrated
against Eisenstein series. The correct 
Eisenstein series to subtract becomes visible from the dominant part 
of the constant term of the Poincar\'e series (below).

Write the Poincar\'e series as 
$$ \Pe(g) \;\;
= \sum_{\gam\in M_k\\G_k} \ph(\gam g) \;\,
= \sum_{\gam\in P_k\\G_k}\; \sum_{\beta\in N_k} \ph(\beta\gam g) 
$$ 
By Poisson summation 
$$ 
 \Pe(g) \;\,
= \sum_{\gam\in P_k\\G_k} \;\, \sum_{\psi \in (N_k\\N_\A)\,\widehat{\;}} \;
\phat_{\gam g}(\psi) \tag 4.4
$$ 
where, $\ph_g(n) = \ph(ng)$, and $\widehat{\varphi}$ is Fourier 
transform along $N_\A$. 
The trivial--$\psi$ (that is, with $\psi=1$) Fourier term 
$$ 
\sum_{\gam\in P_k\\G_k} \phat_{\gam g}(1) \tag 4.5
$$ 
is an Eisenstein series, since the function 
$$ g \to \phat_g(1) 
= \int_{N_\A} \ph(ng)\,dn 
$$ 
is left $M_\A$--equivariant by the character $\delta\chi_0$, and left 
$N_\A$--invariant.

For $\xi\in M_{k},$ 
$$\align  \phat_{\xi g}(\psi)\; 
& =\; \int_{N_\A} \psibar(n)\, \ph(n \xi g)\, dn \\ 
&= \; \int_{N_\A} \psibar(n)\, \ph(\xi \cdot \xi^{-1} n \xi\cdot g)\, dn \;
=\; \int_{N_\A} \psibar(\xi n\xi^{-1})\, \ph(n\cdot g)\, dn \;
=\; \phat_{g}(\psi^\xi) \endalign
$$ 
where $\psi^\xi(n) = \psi(\xi n \xi^{-1})$, by replacing $n$ by $\xi 
n\xi^{-1}$, using the left $M_k$--invariance of $\ph$, and invoking the 
product formula to see that the change-of-measure is trivial. Since 
this action of $Z_k\\M_k$ is transitive on non-trivial characters on 
$N_k\\N_\A$, for a fixed choice of non-trivial character $\psi$, the 
sum over non-trivial characters can be rewritten as a more 
familiar sort of Poincar\'e series 
$$\align 
\sum_{\gam\in P_k\backslash G_k} &\;\,  \sum_{\psi' \in (N_k \backslash N_\A)\,\widehat{\;}} 
\phat_{\gam g}(\psi') \;\,
= 
\sum_{\gam\in P_k \backslash G_k} \;\, \sum_{\xi \in Z_k \backslash M_k} 
\phat_{\gam g}(\psi^\xi) \\ 
&= 
\sum_{\gam\in P_k \backslash G_k} \;\, \sum_{\xi \in Z_k \backslash M_k} 
\phat_{\xi\gam g}(\psi) \;\;\,
= 
\sum_{\gam\in Z_k N_k \backslash G_k} \; 
\phat_{\gam g}(\psi) \endalign
$$ 
Denote this version of the original Poincar\'e series, with the 
Eisenstein series subtracted, by 
$$ 
\Pepsi(g) 
\;\;= 
\sum_{\gam\in Z_k N_k\\G_k} \; \phat_{\gam g}(\psi) 
\;=\; \Pe(g) \;\;- \sum_{\gam\in P_k\\G_k} \phat_{\gam g}(1)\tag 4.6 
$$

\vskip 10pt \noindent {\bf Remark:} With (4.6), the square integrability part of the Poincar\'e series in 
Theorem 2.7 can be precisely formulated as follows. For $\ph$ admissible, the {\it modified} 
Poincar\'e series $\Pepsi(g)$ is in $L^2(Z_\A G_k\\G_\A)$.

\vskip 10pt
Now we describe the continuous part of the spectral decomposition. At
every place $\nu$, let $\eta_\nu$ be the spherical vector in the
(non-normalized) principal series $\text{Ind}_{P_\nu}^{G_\nu}\chi_\nu$,
normalized by $\eta_\nu(1)=1$. Take
$\eta=\bigotimes_{\nu\le \infty}\eta_\nu$. The corresponding
Eisenstein series is  
$$
E_\chi(g) \;\;= \sum_{\gamma\in P_k\\G_k} \eta(\gamma g)
$$
For any left $Z_\A G_k$--invariant and right $K_\A$--invariant
square-integrable $F$ on $G_\A,$ write
$$
\< F, E_\chi \> = \int_{Z_\A G_k\\ G_\A} F(g)\,\overline{E_\chi(g)}\,dg
$$
With suitable normalization of measures,
$$
\hbox{continuous spectrum part of $F$} 
\;\;= 
\int\limits_{\Re(\chi)={1\over 2}} \< F, E_\chi \> \; E_\chi\, d\chi
$$ 
Explicitly, let 
$$ 
\kappa = \text{meas}(\Bbb{J}^{1}/k^{\times}) \tag 4.7
$$ 
the measure on $\Bbb{J}^{1}/k^{\times}$ being the image of the measure $\gamma$ 
on $\Bbb{J}$ defined in \cite{W2}, page 128. 
It is well-known (see \cite{W2}, page 129, Corollary) 
that the residue of the Dedekind zeta-function of $k$ at $s=1$ is 
$$ 
\underset{s = 1}\to{\text{Res}}\; \zeta_{k}(s) 
\;=\; 
{\frac{\kappa}{|D_{k}|^{\frac{1}{2}}}}
$$ 
where $D_{k}$ denotes the discriminant of $k.$ Then, 
the continuous part of $F$ can be written as    
$$
\hbox{continuous spectrum part of $F$} 
\;=\;
{1\over 4\pi i\kappa} 
\sum_\chi \, \int\limits_{\Re(s)=\frac{1}{2}} 
\bra F, E_{s,\chi} \ket 
\cdot E_{s,\chi}\;ds
$$ 
where the sum is over all {\it absolutely} unramified characters $\chi \in \Chat_{0}.$ 
Here $E_{s,\chi}$ stands for $E_{\chi\,|\cdot|^{s}}$ defined above. In general, this formula 
requires isometric extensions to $L^2$ of integral
formulas that converge literally only on a smaller dense subspace
(pseudo-Eisenstein series). However, in our situation, since $\Pepsi(g)$ has sufficient decay, 
its integrals against Eisenstein series (with parameter in a bounded vertical strip containing the 
critical line) converge absolutely. Furthermore, as $\Pepsi$ is sufficiently smooth, 
with derivatives of sufficient decay, its continuous part of the spectral 
decomposition also converges. Then by Theorem 2.7 (see also the above remark), 
Proposition 4.1 and (4.6), with respect to an orthonormal basis $\{F\}$ of 
everywhere locally spherical cuspforms, we have the spectral 
decomposition\footnote{There is no residual contribution to the spectral decomposition of
$\Pepsi(g),$ as can be easily verified.} 
\vskip5pt 
\vbox{ 
$$ 
\Pe 
\;=\;
\left(\int_{N_\infty} \!\!\ph_\infty \right) 
\cdot 
E_{v + 1} 
\;+\;
\sum_{F} \;
\left(\int_{Z_{\infty} \backslash G_{\infty}}\; \ph_{\infty} \cdot \overline{W}_{F,\, \infty} \right)
\cdot 
L(v + {\scriptstyle \frac{1}{2}},\overline{F})\cdot F  \tag 4.8 
$$
$$\hskip62pt \;\;+ \; 
{1\over 4\pi i\kappa} 
\sum_\chi \, \int\limits_{\Re(s)=\frac{1}{2}} 
\bra \Pepsi, E_{s,\chi} \ket 
\cdot E_{s,\chi}\;ds 
$$
}
\noindent where we set $E_{s}:= E_{s, 1}.$ To compute the inner product 
$\bra \Pepsi, E_{s,\chi} \ket$ in the continuous part, first consider an Eisenstein series  
$$ 
E(g)
\;\;= \sum_{\gam\in P_k\\G_k} \eta(\gam g) 
$$ 
for $\eta$ left $P_k$--invariant, left $M_\A$--equivariant and left 
$N_\A$--invariant. The Fourier expansion of this Eisenstein series is 
$$ 
E(g) 
\;\;\,= \sum_{\psi' \in (N_k\\N_\A)\,\widehat{\;}}\;\; 
\int_{N_k\\N_\A} \overline{\psi'}(n) \, E(ng)\,dn 
$$ 
For a fixed non-trivial character $\psi,$ the $\psi^{\text{th}}$ Fourier term is 
$$ 
\int_{N_k\\N_\A} \overline{\psi}(n) \, E(ng)\,dn \,
= \int_{N_k\\N_\A} \psibar(n) \, \sum_{\gam\in P_k\\G_k} \eta(\gam ng)\,dn 
$$ 
$$= \sum_{w\in P_k\\G_k/N_k} \;
\int_{(N_k \, \cap \,  w^{-1}P_k w)\\N_\A} \; \psibar(n) \, \eta(wng)\,dn 
$$ 
$$ 
= 
\int_{N_k\\N_\A} \psibar(n) \, \eta(ng)\,dn \,
\;+\; \int_{N_\A} \psibar(n) \, \eta(w_{\circ} n g)\,dn 
$$
$$= \,
0 \;
+\; \int_{N_\A} \psibar(n) \, \eta(w_{\circ} n g)\,dn
\hskip30pt
\hbox{\Bigg(where $w_{\circ} = \pmatrix 0&1\cr 1 & 0 \endpmatrix$ \Bigg)}
$$
because $\psi$ is non-trivial and $\eta$ is left 
$N_\A$--invariant. Denote the 
$\psi{^{\text{th}}}$ Fourier term by 
$$W^{E}(g) \, =\,  W_{\eta,\,\psi}^{E}(g) \,
= \int_{N_\A} \overline{\psi}(n) \, \eta(w_{\circ} n g)\,dn \tag 4.9 
$$

We have the following 

\vskip 12pt
\proclaim{Proposition 4.10} Fix $s\in \Bbb{C}$ such that $\Re(s) > 1,$
and suppose 
$\ph_\infty \in \Omega_\ell$ with $\Re(v),\, \ell$ sufficiently large. 
Then,   
$$  
\bra \text{P\'e}\,^{*}, E_{s,\chi} \ket \,=\,\overline{\chi}(\diff)   
\left(\int_{Z_{\infty} \backslash G_{\infty}} \ph_{\infty} \cdot 
\overline{W}_{s,\, \chi,\, \infty}^{E} \right)
{L(v + \bar{s}, \overline{\chi}) 
\cdot 
L(v + 1 - \bar{s}, \chi) 
\over 
L(2\bar{s},\overline{\chi}^2)} 
\cdot 
|\diff|^{-(v - \bar{s} + 1/2)} 
$$ 
where $\diff$ denotes a differental idele (see \cite{W2}, page 113, Definition 4) with component $1$ 
at archimedean places.
\endproclaim

\vskip 10pt 
\noindent{\bf Proof:} Fix a non-trivial character $\psi$ on $N_{k}\backslash N_{\A}$. 
When $\Re(v)$ and $\ell$ are both large, the modified Poincar\'e series $\Pepsi(g)$ has 
sufficient (polynomial) decay, and therefore,
we can unwind it to obtain (see (4.6)): 
$$ 
\int_{Z_\A G_k\\G_\A} \Pepsi(g)\, \overline{E}_{s,\chi}(g)\, dg \,
=  \int_{Z_\A N_\A\\G_\A} \int_{N_k\\N_\A} \phat_{ng}(\psi)\, 
\overline{E}_{s,\chi}(ng)\,dn\, dg \tag 4.11 
$$ 
$$ 
= \, \int_{Z_\A N_\A\\G_\A} \phat_g(\psi) 
\int_{N_k\\N_\A} \psi(n) \overline{E}_{s,\chi}(ng)\,dn\, dg 
\, =  \int_{Z_\A N_\A\\G_\A} \phat_g(\psi)\, \overline{W}_{s,\chi}^{E}(g)\, dg  
$$ 
Since 
$$
\phat_{g}(\psi)\; 
 =\; \int_{N_\A} \psibar(n)\, \ph(n g)\, dn
$$ 
the last integral in (4.11) is 
$$\align 
\int_{Z_\A N_\A \backslash G_\A} \int_{N_\A} \psibar(n)\, \ph(n g)\, 
\overline{W}_{s,\chi}^{E}(g)\, dn\, dg \;
&=\; \int_{Z_\A N_\A\backslash G_\A} \int_{N_\A} 
\ph(n g)\, \overline{W}_{s,\chi}^{E}(ng)\, dn\, dg\\
& = \int_{Z_\A \backslash G_\A}  \ph(g)\, \overline{W}_{s,\chi}^{E}(g)\, dg \tag 4.12
\endalign$$

The Whittaker function of the Eisenstein series does factor over 
primes, into local factors depending only upon the local 
data at $\nu$ 
$$ 
W^{E}_{s,\chi} \,
= \,\bigotimes_{\nu} W^{E}_{s,\chi, \nu} 
$$ 
Thus, by (4.11) and (4.12),
$$ 
\bra \Pepsi, E_{s,\chi}\ket \,
= \,
\left(\int_{Z_{\infty}\\G_{\infty}} \ph_{\infty}\cdot \overline{W}^E_{s,\chi, \infty}\right) 
\cdot \prod_{\nu < \infty} \int_{Z_{\nu}\\G_{\nu}} \ph_{\nu}(g_{\nu})\,\overline{W}^E_{s,\chi, \nu}(g_{\nu})\,dg_{\nu} 
$$

At finite $\nu,$ using an Iwasawa decomposition and the vanishing of 
$\ph_{\nu}$ off $M_{\nu} K_{\nu}$ (see (3.2)), as in the integration against cuspforms, the 
local factor is 
$$ 
\int_{k_{\nu}^{\times}} |a|_{\nu}^{v}\,\overline{W}^{E}_{s,\chi, \nu}\pmatrix a&\cr &1\endpmatrix\,da 
$$ 
However, for Eisenstein series, the natural normalization of the 
Whittaker functions differs from that used for cuspforms, instead 
presenting the local Whittaker functions as images 
under intertwining operators. Specifically, define the {\it 
normalized spherical vector} for data $s,\, \chi_{\nu}$ 
$$ 
\eta_{\nu}(p k) = |a/d|_{\nu}^{s}\cdot \chi_{\nu}(a/d) 
\qquad \left(\text{for $p=\pmatrix a&*\cr &d \endpmatrix \in P_{\nu}$ and $k \in K_{\nu}$} \right) 
$$ 
The corresponding spherical local Whittaker function 
for Eisenstein series is the integral 
\footnote 
{This integral only converges nicely for $\Re(s)\gg 0,$ but admits a 
meromorphic continuation in $s$ by various means. For example, the 
algebraic form of Bernstein's {\it continuation 
principle} applies, since the dimension of the space of intertwining 
operators from the principal series to the Whittaker space is 
one-dimensional. 
} (see (4.9))
$$ 
W^{E}_{s, \chi, \nu}(g) 
= 
\int_{N_{\nu}} \overline{\psi}_{\nu}(n)\,\eta_{\nu}(w_{\circ} ng)\,dn 
$$ 
The Mellin transform of 
the Eisenstein-series normalization $W^{E}_{s, \chi, \nu}$ compares 
to the Mellin transform of the usual normalization as follows. Let 
$\diff_{\nu}\in k_{\nu}^\times$ be such that 
$$ 
(\o_{\nu}^{*})^{-1} = \diff_{\nu}\cdot \o_{\nu} 
$$ 
Let $\diff$ be the idele with $\nu^{\text{th}}$ component $\diff_{\nu}$ at finite 
places $\nu$ and component $1$ at archimedean places. 
Then for finite $\nu$ the $\nu^{\text{th}}$ local integral is (see Appendix 2 for details), 
$$ 
\int_{k_{\nu}^{\times}}\, |a|_{\nu}^{v}\,\overline{W}^{E}_{s, \chi, \nu} \pmatrix a&0\cr 0&1\endpmatrix\,da 
\,=\, 
|\diff_{_\nu}|_{\nu}^{1/2} 
\cdot
{ 
L_{\nu}(v + \bar{s}, \overline{\chi}_{\nu}) 
\cdot 
L_{\nu}(v + 1 - \bar{s}, \chi_{\nu}) 
\over 
L_{\nu}(2\bar{s},\overline{\chi}_{\nu}^2)  
}
\cdot 
|\diff_{_\nu}|_{\nu}^{-(v + 1 - \bar{s})}\,\overline{\chi}_{\nu}(\diff_{_\nu}) 
$$ 
and the proposition follows. \qquad 
\qed

Accordingly, the spectral decomposition (4.8) is  
 $$\align 
&\Pe 
\;=\;
\left(\int_{N_\infty} \!\!\ph_\infty \right) 
\cdot 
E_{v + 1} 
\;+\; 
\sum_{F} \;
\left(\int_{Z_{\infty} \backslash G_{\infty}}\; \ph_{\infty} \cdot \overline{W}_{F,\, \infty} \right)
\cdot 
L(v + {\scriptstyle \frac{1}{2}},\overline{F})\cdot F  \tag 4.13\\
& +  
\sum_{\chi} \, \frac{ \overline{\chi}(\diff)}{4\pi i\kappa} \int\limits_{\Re(s)=\frac{1}{2}}   
\left(\int_{Z_{\infty} \backslash G_{\infty}} \ph_{\infty} \cdot 
W_{1 - s,\, \overline{\chi},\, \infty}^{E} \right)
{L(v + 1 - s, \overline{\chi}) 
\cdot 
L(v + s, \chi) 
\over 
L(2 - 2s,\overline{\chi}^2)} 
\, 
|\diff|^{-(v + s - 1/2)} 
 \cdot E_{s,\chi}\;ds 
\endalign$$
where we replaced $\bar{s}$ by $1 - s,$ for $\Re(s) = \frac{1}{2},$ to maintain holomorphy 
of the integrand. The archimedean-place Whittaker functions can be expressed in terms of 
the usual $K$--Bessel function as follows. Let 
$$ 
\eta_{\nu}(n m k) = \eta_{s,\, \nu}(n m k) = |a/d|^{s}_{\nu} 
\qquad \left(\text{for $n\in N_{\nu},$ $m=\pmatrix a&\cr &d \endpmatrix \in M_{\nu},$ 
$k \in K_{\nu}$}\right) 
$$ 
The normalization of the Whittaker function is 
$$ 
W^{E}_{s,\, \nu}(g) = \int_{N_{\nu}} \psibar_{\nu}(n)\,\eta_{\nu}(w_{\circ} n g)\,dn 
\qquad \left(\text{for $\Re(s)\gg 0$ and fixed non-trivial $\psi$}\right) 
$$ 
Then, for $\nu$ archimedean and fixed non-trivial character $\psi_{0, \, \nu}$ 
on $k_{\nu}$   
$$ 
W^{E}_{s,\, \nu}\pmatrix a&\cr &1\endpmatrix 
= \int_{k_{\nu}} \psibar_{0,\, \nu}(x)\, 
\left| {a\over aa^{\iota} + xx^{\iota}}\right|_{\nu}^{s} 
\,dx 
= |a|_{\nu}^{1 - s} \int_{k_{\nu}} \psibar_{0, \, \nu}(ax)\, 
{1\over |1 + xx^{\iota}|_{\nu}^{s}} 
\,dx 
$$ 
by replacing $x$ by $ax,$ where $\iota$ is the complex conjugation for 
$\nu \isom \C$ and the identity map for $\nu \isom \R.$ The usual 
computation shows that 
$$ 
W^{E}_{s,\, \R}\pmatrix a&\cr &1\endpmatrix  
= 
{|a|^{1/2} \over \pi^{-s}\Gamma(s)} 
\int_{0}^{\infty} e^{- \pi (t + {1\over t})|a|}\,t^{s - \frac{1}{2}}\,{dt\over t} 
= 
{2\,|a|^{1/2} K_{s  - 1/2}(2 \pi |a|)\over \pi^{-s}\Gamma(s)} 
$$ 
and, similarly, 
\footnote
{The appropriate measure is {\it double} the 
usual in this case. 
} 
$$ 
W^{E}_{s,\, \C}\pmatrix a&\cr &1 \endpmatrix  
= 
{|a| \over (2\pi)^{-2s}\Gamma(2s)} 
\int_{0}^{\infty} e^{- 2 \pi (t + {1\over t})|a|}\,t^{2s - 1}\,{dt\over t} 
= 
{2\, |a| K_{2 s  - 1}(4 \pi |a|)\over (2\pi)^{-2s}\Gamma(2s)}  
$$

To simplify the integral over $Z_{\infty}\backslash G_{\infty}$ in the continuous part 
of (4.13), let 
$$ 
\Phi_{\nu}(x) = \ph_{\nu}\pmatrix 1&x\cr 0 &1 \endpmatrix \qquad (\text{for $\nu$ archimedean})
$$ 
Using the right $K_{\nu}$--invariance and an Iwasawa decomposition, 
$$\align 
\int_{Z_{\nu}\backslash G_{\nu}} \ph_{\nu} 
\cdot 
W^E_{s,\, \nu} 
&= 
\int_{k_{\nu}^{\times}}\int_{k_{\nu}} 
|a|_{\nu}^{v}\, \Phi_{\nu}(x)\, W^{E}_{s,\, \nu}\pmatrix a&\cr &1\endpmatrix 
\psi_{0, \, \nu}(ax) 
\, 
dx\,da \tag 4.14\\  
& = 
\int_{k_{\nu}^{\times}} |a|_{\nu}^{v}\,
\Phihat_{\nu}(a)\, W^{E}_{s,\, \nu}\pmatrix a&\cr &1\endpmatrix da
\endalign$$ For $\chi \in \Chat_{0}$ absolutely unramified, we have 
$W_{s,\, \chi,\, \nu}^{E} = W_{s + it_{\nu},\, \nu}^{E},$ where $t_{\nu} \in \R$ is 
the parameter of the local component $\chi_{\nu}$ of $\chi.$ Then all archimedean 
integrals in the continuous part of the spectral decomposition of $\Pepsi$ 
are given by (4.14) with $s$ replaced by $1 - s -  it_{\nu}.$

In particular, if $\ph_{\nu}$ is specialized to (3.11) we have 
$$\int_{Z_{\nu}\backslash G_{\nu}} \ph_{\nu} 
\cdot 
W^E_{s,\, \nu} \, 
=\, \cases {\displaystyle \frac{\Cal{G}_{\nu}(s; v,
w)}{\pi^{-s}\Gamma(s)}} & \text{if $\nu \isom \R$}\\
 {\displaystyle \frac{\Cal{G}_{\nu}(s; v,
w)}{(2\pi)^{-2s - 1}\Gamma(2s)}} & \text{if $\nu \isom \C$}
\endcases \tag 4.15$$ 
where $\Cal{G}_{\nu}(s; v, w)$ is given in (4.2) and (4.3). Furthermore, with 
these choices of $\ph_{\nu},$ 
$$\int_{N_{\nu}} \!\!\ph_{\nu}\, =\, 
 \cases \sqrt{\pi}\, {\displaystyle \frac{\Gamma(\frac{w - 1}{2})}{ \Gamma(\frac{w}{2})}} 
 & \text{if $\nu \isom \R$}\\
2\pi\, (w - 1)^{-1} & \text{if $\nu \isom \C$}
\endcases \tag 4.16$$

As usual, let $r_1$ and $r_2$ denote the number of
real and complex embeddings of $k,$ respectively. Following \cite{DG1}, 
Proposition 5.10, we now prove

\vskip 12pt
\proclaim{Theorem 4.17} Assume $\ph$ is defined by (3.11). The
Poincar\'e series $\Pe(g)$ has meromorphic continuation to a region
in $\Bbb{C}^2$ containing $v = 0,$ $w = 1.$ As a function of
$w,$ for $v =0,$ it is holomorphic in the half-plane $\Re(w) > 11/18$, 
except for $w = 1$ where it has a
pole of order $r_1 + r_2 + 1.$ 
\endproclaim

\vskip 10pt 
\noindent{\bf Proof:} Let $\Pepsi_{\text{cusp}}$ and $\Pepsi_{\text{cont}}$ be, respectively, 
the discrete and continuous parts of $\Pepsi.$ Then the spectral decomposition (4.13) is 
$$ 
\Pe 
\;=\;
R(w) 
\cdot 
E_{v + 1} \;+\; \Pepsi_{\text{cusp}} \;+\; \Pepsi_{\text{cont}} \qquad  
\left(\text{where $
R(w) 
=
\int_{N_\infty} \!\!\ph_\infty
$} \right)
$$ 
the integral being computed by (4.16). As mentioned before (see the discussion 
after the proof of Proposition 4.1), 
the series giving $\Pepsi_{\text{cusp}}$ 
$$
\sum_{F} \; \bar{\rho}_{_F}\,  \Cal{G}_{_{F_{_\infty}}}(v,
w) \, L(v + {\scriptstyle \frac{1}{2}},\overline{F})\cdot F  
$$ 
converges absolutely for 
$(v, w)\in \C^{2},$ apart from the poles of 
$\Cal{G}_{\nu}({\scriptstyle \frac{1}{2}} + i\bar{\mu}_{_{F, \nu}}; v,
w).$ The fact that $\Pepsi_{\text{cusp}}$ equals the above spectral sum 
is justified by the square integrability of $\Pepsi$ for $\Re(w) > 1$ 
and large $\Re(v)$ (see Theorem 2.7, (4.6) and Appendix 1). Furthermore, using (4.2), (4.3) 
and the Kim-Shahidi bound of the local parameters $|\Re(i\mu_{_{f, \nu}})| < 1/9$ 
(see \cite{K}, \cite{KS}), 
the cuspidal part $\Pepsi_{\text{cusp}}$ is 
holomorphic for $v = 0$ and $\Re(w) > 11/18$.

To deal with the continuous part, first note that by (4.15) the expression under 
the vertical line integral in (4.13) has enough decay in the parameters to ensure 
the absolute convergence of the integral and sum over $\chi.$ Also, note that 
$\Pepsi_{\text{cont}}$ is holomorphic for $\Re(v) > \frac{1}{2}$ and $\Re(w) > 1.$ 
Aiming to analytically continue to $v = 0,$ first take $\Re(v)= 1/2 + \varepsilon,$ 
and move the line of integration from 
$\sigma = 1/2$ to $\sigma = 1/2 - 2\varepsilon.$ This picks up the 
residue of the integrand corresponding to $\chi$ trivial due to the 
pole of $\zeta_{k}(v + s)$ at $v + s=1,$ that is, at $s=1- v.$ Its contribution is 
$$  
{1\over 2}\, 
Q(v; v, w) 
\;\cdot \; 
|\diff|^{1/2} 
\cdot 
|\diff|^{-1/2} 
\cdot 
E_{1 - v} 
\,=\, 
{1\over 2}\,
Q(v; v, w)   
\cdot 
E_{1 - v} 
$$ 
where 
$$
Q(s; v, w)
\,=\, 
\int_{Z_\infty\\G_\infty}\!\! \ph_\infty\cdot W^{E}_{s,\,  \infty}
$$ 
stands for the ratio of products of gamma functions computed by (4.15). 
This expression of $\Pepsi_{\text{cont}}$ is holomorphic in $v$ in the strip 
$$ 
\frac{1}{2} - \varepsilon \le \Re(v) \le \frac{1}{2} + \varepsilon
$$

Now, take $v$ with $\Re(v) = 1/2 - \varepsilon,$ and then move the vertical integral 
from $\sigma = 1/2 - 2\varepsilon$ back to $\sigma = 1/2.$ 
This picks up $(-1)$ times the residue at the pole of 
$\zeta_{k}(v + 1 - s)$ at $1,$ that is, at $s = v,$ with another sign due to 
the sign of $s$ inside this zeta function. Thus, we pick up the 
residue  
$$  
{1\over 2}\,
Q(1 - v; v, w)   
\cdot 
{ 
\zeta_{k}(2 v) 
\over 
\zeta_{k}(2 - 2v) 
} 
\cdot 
|\diff|^{-2v+1} 
\cdot 
E_{v}\,
 = \,
{1\over 2}\, 
Q(1 - v; v, w)   
\cdot 
{ 
\zeta_{\infty}(2 - 2 v) 
\over 
\zeta_{\infty}(2v) 
} 
\cdot 
E_{1 - v}
$$ 
where the last identity was obtained from the functional equation of the Eisenstein series 
$E_{v}.$ Since $\Cal{G}_{\nu}(s; v, w)$ defined in (4.2) and (4.3) is invariant under $s \to 1 - s,$ 
it follows by (4.15) that the above residues are equal. Note that the part of $\Pepsi_{\text{cont}}$ 
corresponding to the vertical line integral and the sum over $\chi$ is now holomorphic in a 
region of $\Bbb{C}^{2}$ containing $v = 0,$ $w = 1.$ In particular, for $v = 0,$ this part of the 
continuous spectrum is holomorphic in the half-plane $\Re(w) > 1/2.$

On the other hand, by direct computation, the {\it apparent} pole of $R(w)E_{v + 1}$ 
at $v = 0$ (independent of $w$) cancels with the corresponding pole of $Q(v; v, w)E_{1 - v}$. 
To establish that the order of the pole at $w = 1$, when $v = 0$, is $r_{1} + r_{2} + 1$, 
consider the most relevant terms (recall (4.15), (4.16)) in the Laurent expansions of 
$R(w)E_{v + 1}$ and $Q(v; v, w)E_{1 - v}$. Putting them together, we obtain an expression  
$$\frac{1}{v}\cdot \left[\frac{c_{1}}{(w - 1)^{r_{1} + r_{2}}} \;
-\; \frac{c_{2}}{(2v + w - 1)^{r_{1} + r_{2}}} \right]$$ 
for some constants $c_{1},\, c_{2}.$ As there is no pole at $v = 0,$ we have 
$c_{1} = c_{2}.$ Canceling the factor $1/v,$ and then setting $v = 0,$ the assertion 
follows. 

This completes the proof. \qquad 
\qed

\vskip 20pt \noindent {\bf \S 5. Asymptotic formula} 
\vskip 10pt

Let $k$ be a number field with $r_1$ real embeddings and $2r_2$
complex embeddings.  Assume that $\ph$ is specialized to (3.11). By
Theorem 3.12, for
$\Re(v)$ and $\Re(w)$ sufficiently large, the integral $I(\chi_0) =
I(v, w)$ defined by (3.6) is
$$I(v, w)  \; =   \sum_{\chi \in \Chat_{0, S}}\; \frac{1}{2\pi i}\int\limits_{_{\Re(s) =\sigma}} L( \chi^{-1}| \cdot |^{v + 1 - s},\,f_1) 
\cdot 
L(\chi |\cdot|^{s},\,\bar{f}_{2})\, 
\Cal{K}_{\infty}(s,\, v,\,w,\,  \chi) \, ds\tag 5.1$$ where $\Cal{K}_{\infty}(s,\, v,\,w,\,  \chi)$ is given by (3.9) 
and (3.10), and where the sum is over $\chi \in \Chat_{0}$ unramified
 outside $S$ and with bounded ramification, depending only on $f_1$
 and $f_2$.

 By Theorem 4.17, it follows that $I(v, w)$ admits meromorphic continuation to a region in $\Bbb{C}^2$ containing the point $v = 0$, $w = 1$. In particular, if $f_{1} = f_{2} = \bar{f}$, then $I(0, w)$ is holomorphic for $\Re(w) > 11/18$, except for $w = 1$ where it has a pole of order $r_1 + r_2 + 1$.

We will shift the line of integration to 
$\Re(s) = \frac{1}{2}$ in (5.1) and set $v = 0$. To do so, we need
{\it some} analytic continuation and reasonable decay in
$|\Im(s)|$ for the kernel function $\Cal{K}_{\infty}(s,\, v,\,w,\,
\chi)$. In fact, it is desirable for applications to have {\it
precise} asymptotic formulae as the parameters $s,\,v,\,  w,\, \chi$
vary. By the decomposition (3.10), the analysis of the
kernel $\Cal{K}_{\infty}(s,\, v,\,w,\,  \chi)$ reduces to
the corresponding analysis of the local component $\Cal{K}_{\nu}(s,\,
v,\,w,\,  \chi_{\nu})$, for $\nu | \infty$. When $\nu$ is complex, one
can use the asymptotic formula already established in \cite{DG2},
Theorem 6.2. For coherence, we include a
simple computation matching, as it should, the local integral (3.9),
for $\nu$ complex, with the integral (4.15) in \cite{DG2}.

Fix a complex place $\nu| \infty$. An irreducible unitary representation of $GL_{2}(\C)$ 
always contains a spherical vector. Therefore, since we are interested only in the 
finite-prime part of the $L$--function associated to a cuspidal 
representation, we may as well suppose that $f_1$ and $f_2$ are spherical at all complex 
places. Also, recall that any character $\chi_{\nu}$ of 
$Z_{\nu}\backslash M_{\nu}\isom \Bbb C^{\times}$ has the form 
$$
\chi_{\nu}(m_{\nu}) \,=\, |z_{\nu}|_{_\Bbb{C}}^{^{\frac{\ell_{\nu}}{2} + i t_{\nu}}} z_{\nu}^{-\ell_{\nu}} \qquad \;\;\;\left(m_{\nu} = \pmatrix
z_{\nu}& 0\\ 0 & 1\endpmatrix,\; t_{\nu}\in \Bbb{R},\; \ell_\nu\in \Bbb{Z}\right)
$$ 
Then, the local integral (3.9) at $\nu$ is 
$$\align &\Cal{K}_{\nu}(s,\, v,\,w,\,  \chi_{\nu})\,
=\, \int\limits_{0}^{\infty}\, \int\limits_{0}^{\infty}\,  \int\limits_{\Bbb{C}}\;
\int\limits_{-\pi}^{\pi} \,\int\limits_{-\pi}^{\pi}\; (|x|^2 + 1)^{^{-w}}\, e^{2\pi i \cdot \text{Tr}_{\Bbb{C}/{\Bbb{R}}}(a_{_1} x e^{i\theta_1} - a_{_2} x e^{i\theta_2})}\\
& \hskip 58pt\cdot a_{1}^{2v + 1 - 2s - 2it_{\nu}} K_{2 i\mu_{_1}}(4 \pi a_1)\, a_{2}^{2s + 2it_{\nu}  - 1} K_{2 i \bar{\mu}_{_2}}(4 \pi a_2)\, e^{i\ell_\nu \theta_1} e^{-i\ell_\nu \theta_2}
\; d\theta_1  d\theta_2  \,  dx \, da_{1} da_{2}\endalign$$
Replacing $x$ by $x/a_{1}$, we obtain 
$$\align &\Cal{K}_{\nu}(s,\, v,\,w,\,  \chi_{\nu}) \,
= \, \int\limits_{0}^{\infty} \, \int\limits_{0}^{\infty} \, \int\limits_{\Bbb{C}}\;
\int\limits_{-\pi}^{\pi} \, \int\limits_{-\pi}^{\pi}\; \left(\frac{a_1}{\sqrt{|x|^2 + a_{1} ^2}}\right)^{2w} e^{2\pi i \cdot \text{Tr}_{\Bbb{C}/{\Bbb{R}}}\big(x e^{i\theta_1} - \frac{a_{_2}}{a_{_1}} x e^{i \theta_2}\big)}\\
& \hskip 58pt\cdot a_{1}^{2v - 1 - 2s - 2it_{\nu}} K_{2 i \mu_{_1}}(4 \pi a_1)\, a_{2}^{2s + 2it_{\nu} - 1} 
K_{2 i \bar{\mu}_{_2}}(4 \pi a_2)\, e^{i\ell_\nu \theta_1} e^{-i\ell_\nu \theta_2}
\; d\theta_1  d\theta_2  \,  dx \, da_{1} da_{2}\endalign$$
If we further substitute 
$$a_1 = r\cos \phi \quad x_1 = r \sin \phi \cos \theta \quad
x_2 = r \sin \phi \sin \theta \quad a_2 = u \cos \phi$$ with 
$0\le \phi \le \frac{\pi}{2}$ and $0\le \theta \le 2\pi$, then 
$$\align &\Cal{K}_{\nu}(s,\, v,\,w,\,  \chi_{\nu}) \,
=\, \int\limits_{0}^{\infty}\, \int\limits_{0}^{\infty} \, \int\limits_{0}^{\frac{\pi}{2}} \, \int\limits_{0}^{2\pi}\,
\int\limits_{-\pi}^{\pi}\, \int\limits_{-\pi}^{\pi}\; (\cos \phi)^{2w + 2v  - 1}\, e^{2\pi i \cdot \text{Tr}_{\Bbb{C}/{\Bbb{R}}}(r\sin \phi  \cdot  e^{i(\theta + \theta_1)} -  u \sin \phi \cdot  e^{i(\theta + \theta_2)})}\\
& \hskip 3pt\cdot r^{2v + 1 - 2s - 2it_{\nu}} K_{2 i\mu_{_1}}(4 \pi r \cos \phi )\, u^{2s + 2it_{\nu} - 1} 
K_{2 i \bar{\mu}_{_2}}(4 \pi u \cos \phi)\, e^{i\ell_\nu \theta_1} e^{-i\ell_\nu \theta_2} \sin \phi
\; d\theta_1  d\theta_2  \,  d\theta \, d\phi \, dr \, du\endalign$$
Using the Fourier expansion 
$$e ^{ i t\sin \theta } 
\; =  \; \sum_{k = -\infty} ^{\infty} J_{k}(t)\, e ^{i k
\theta}$$ we obtain 
$$\align \Cal{K}_{\nu}(s,\, v,\,w,\,  \chi_{\nu}) =  (2 \pi)^3 \int\limits_{0}^{\infty}&\int\limits_{0}^{\infty}
\int\limits_{0}^{\frac{\pi}{2}} K_{2 i\mu_{_1}}(4\pi r\cos \phi) K_{2 i \bar{\mu}_{_2}}(4\pi u\cos \phi)
J_{\ell_{_\nu}} (4\pi r\sin \phi) J_{\ell_{_\nu}} (4\pi u\sin \phi)
\\ & \hskip 70pt \cdot u^{2s + 2it_{\nu}} r^{2v + 2 - 2s - 2it_{\nu}} (\cos \phi)^{2w  +  2v  - 1} \sin \phi \;
\frac{d\phi dr du}{ru}\endalign$$
In the notation of \cite{DG2}, equation (4.15), this is essentially 
$\Cal K_{\ell_{_{\nu}}}(2s + 2it_{\nu}, 2v, 2w)$. It follows that $\Cal{K}_{\nu}(s, v, w,  \chi_{\nu})$ 
is analytic in a region $\Cal{D}:$ 
 $\Re(s) = \sigma > \frac{1}{2}  - \varepsilon_{0}$, $\Re(v) > -\varepsilon_{0}$ and $\Re(w) > \frac{3}{4}$, with a fixed (small) $\varepsilon_{0} > 0$, and moreover, we have the asymptotic formula
 $$\align \Cal{K}_{\nu}(s,\, v,\,w,\,  \chi_{\nu})  = \pi^{-2v + 1} A(v, w, \mu_{_{1}},\mu_{_{2}})& \cdot
 \big(1 + \ell_{\nu}^2+ 4(t + t_{\nu})^2 \big)^{-w}\\
&  \cdot \bigg[1 \; + \; \Cal O_{\sigma,\, v,\, w,\, \mu_{_{1}},\, \mu_{_{2}}}\left(
\Big(\sqrt{1 + \ell_{\nu}^2+ 4(t + t_{\nu})^2}\Big)^{-1}\right)\bigg]\tag 5.2\endalign$$ 
where $A(v,w,\mu_{1},\mu_{2})$ is the ratio of products of gamma functions 

$$2^{4w - 4v - 4}\,
\frac{\Gamma(w+v+i\mu_{1}+i \bar{\mu}_{2})
 \Gamma(w+v-i\mu_{1} + i \bar{\mu}_{2}) 
\Gamma(w+v+i\mu_{1} - i \bar{\mu}_{2})
\Gamma(w+v-i\mu_{1} - i \bar{\mu}_{2})}{\Gamma(2w + 2v)}\tag 5.3$$

For $\nu$ real, the corresponding argument (including the integrals that arise from the (anti-) holomorphic discrete series) is even simpler (see
\cite{DG1} and \cite{Zh2}). In this case, the  
asymptotic formula of $\Cal{K}_{\nu}(s, v, w,  \chi_{\nu})$ becomes 
 $$\align \Cal{K}_{\nu}(s,\, v,\,w,\,  \chi_{\nu})  =  B(v, w, \mu_{_{1}},\mu_{_{2}})& \cdot
 \big(1 + |t + t_{\nu}| \big)^{-w}\\
&  \cdot \bigg[1 \; + \; \Cal O_{\sigma,\, v,\, w,\, \mu_{_{1}},\, \mu_{_{2}}}\left(
\big(1 + |t + t_{\nu}|\big)^{-\frac{1}{2}}\right)\bigg]\tag 5.4\endalign$$ where $B(v, w, \mu_{_{1}},\mu_{_{2}})$ is a similar ratio of products of gamma functions. 

It now follows that for $\Re(w)$ sufficiently large,
$$I(0, w)  \; =   \sum_{\chi \in \Chat_{0, S}}\; \frac{1}{2\pi}\int\limits_{-\infty}^\infty L\big( \chi^{-1}| \cdot |^{\frac{1}{2} - it},\,f_1\big) 
\cdot 
L\big(\chi |\cdot|^{\frac{1}{2} + it},\,\bar{f}_{2}\big)\, 
\Cal{K}_{\infty}({\scriptstyle \frac{1}{2}} + it,\, 0,\,w,\,  \chi) \,
dt\tag 5.5$$ 
Since $I(0, w)$ has analytic continuation to $\Re(w) > 11/18$, a mean value result 
can already be established by standard arguments. For instance, assume $f_{1} = f_{2} = \bar{f}$, and 
choose a function $h(w)$ which is holomorphic and with sufficient decay (in $|\Im(w)|$) 
in a suitable vertical strip containing $\Re(w) = 1$. For example, one can choose a suitable product of gamma functions. Consider the integral    
$$\frac{1}{i}\,\int\limits_{_{\Re(w)  = L}} I(0, w)\,h(w)\,T^{w} \, dw\tag 5.6$$ with $L$ a large positive constant. Assuming $h(1) = 1$, we have the asymptotic formula
$$\sum_{\chi \in \Chat_{0, S}}\;\; \int\limits_{-\infty}^{\infty}\; 
|L({\scriptstyle \frac{1}{2}} + it,  f \otimes \chi)|^2\cdot M_{_{\chi,\,T}}(t)\, dt \; \sim \; A\,T\,(\log T)^{^{r_1 + r_2}}\tag 5.7$$ for some computable positive constant $A$, where
$$M_{_{\chi,\,T}}(t) \, = \, \frac{1}{2\pi i}\,\int\limits_{_{\Re(w)  = L}} \Cal{K}_{\infty}({\scriptstyle \frac{1}{2}} + it,\, 0,\,w,\,  \chi) 
\,h(w)\,T^{w} \, dw\tag 5.8$$

For a character $\chi \in \Chat_{0},$ put 
$$\kappa_{_{\chi}}(t)\,= \, \prod_{\nu \isom \Bbb{R}}\,
\left(1 + |t + t_{\nu}|\right) \, \cdot \,
\prod_{\nu \isom \Bbb{C}}\,\left(1 + \ell_{\nu}^2+ 4(t + t_{\nu})^2 \right) 
\quad\;\;  (t\in \Bbb{R}) 
\tag 5.9$$ 
where $it_{\nu}$ and $\ell_{\nu}$ are the parameters of the
local component $\chi_{\nu}$ of $\chi$. Since $\chi$ is trivial on the 
positive reals, 
$$\sum_{\nu| \infty} \,d_{\nu}\, t_{\nu} = 0$$ 
with $d_{\nu} = [k_\nu:\R]$ the local degree. Then, the main
contribution to the asymptotic formula (5.7) comes from terms for
which $\kappa_{_{\chi}}(t)\ll T.$ For applications, it might   
be more convenient to work with a slightly modified function $Z(w)$ defined by   
$$Z(w)  \; =   \sum_{\chi \in \Chat_{0, S}}\;\; \int\limits_{-\infty}^\infty\;  
|L({\scriptstyle \frac{1}{2}} + it, f \otimes \chi)|^{2} \cdot 
\kappa_{_{\chi}}(t)^{-w} \,
dt \tag 5.10$$ 
It is obtained from the function $I(0, w)$ by essentially picking off just the main terms in the asymptotic 
formulae (5.2) and (5.4) of the local components $\Cal{K}_{\nu}(s, 0, w,  \chi_{\nu}).$ Its analytic properties can be transferred (via some technical adjustments) from those of $I(0, w).$ As an 
illustration of this fact, we show that the right-hand side of (5.10) is absolutely convergent for 
$\Re(w) > 1.$ Using the asymptotic formulae (5.2) and (5.4), it
clearly suffices to verify the absolute convergence of the right-hand
side of (5.5), with $f_{1} = f_{2} = \bar{f},$ when $w > 1.$ 

To see the absolute convergence of the defining expression (5.5) for
$I(0,w)$, first note that the triple integral expressing
$\Cal{K}_{\nu}(s,\, v,\,w,\,  \chi_{\nu})$ can be written as    
$$
\Cal{K}_{\nu}({\scriptstyle \frac{1}{2}} + it,\, 0,\,w,\,  \chi_{\nu}) \,=\,  (2 \pi)^3 
\int\limits_{0}^{\frac{\pi}{2}}\, (\cos \phi)^{2w - 1} \sin \phi 
\cdot | V_{\mu_{_{f, \nu}},\, \chi_{_{\nu}}}(t, \phi)|^{2}\,d\phi \qquad (\text{for $\nu \isom \Bbb{C}$})\tag 5.11
$$ 
when $v = 0$ and $\Re(s) = \frac{1}{2}$, where 
$$
V_{\mu_{_{f, \nu}},\, \chi_{_{\nu}}}(t, \phi)\,=\,\int\limits_{0}^{\infty}\, u^{2i(t_{\nu} + t)} K_{2 i\mu_{_{f, \nu}}}(4\pi u\cos \phi)
J_{_{|\ell_{_\nu}|}} (4\pi u\sin \phi)\, du \tag 5.12
$$
Here we also used the well-known identity $J_{-\ell_{\nu}}(z) = (-1)^{\ell_{\nu}}J_{\ell_{\nu}}(z).$ The 
convergence of the last integral is justified by 6.576, integral 3, page 716 in \cite{GR}. For 
$\nu \isom \Bbb{R},$ the local integral (3.9) has a similar form, when $v = 0$ 
and $\Re(s) = 1/2,$ as it can be easily verified by a straightforward computation.

The form of the integral (5.11) allows us to adopt the argument used in the proof of Landau's Lemma 
to our context giving the desired conclusion. We shall follow \cite{C}, proof of Theorem 6, page 115.

Choose a sufficiently large real number $a$ such that the right-hand side of (5.5) is convergent at 
$w = a.$ Since $I(0, w)$ is holomorphic for $\Re(w) > 1,$ its Taylor series 
$$\align
&\sum_{j = 0}^{\infty}\, {(w - a)^{j}\over j!}\, I^{(j)}(0, a) \tag 5.13\\
&=\, 
\frac{1}{2\pi}\sum_{j = 0}^{\infty}\,\, {(w - a)^{j}\over j!}
\sum_{\chi \in \Chat_{0, S}}\;\, \int\limits_{-\infty}^\infty\, 
|L({\scriptstyle \frac{1}{2}} + it, f \otimes \chi)|^{2}\cdot  
\Cal{K}_{\infty}^{(j)}({\scriptstyle \frac{1}{2}} + it,  0, a,  \chi) \,
dt 
\endalign 
$$ 
has radius of convergence $a - 1.$ Using the structure of (5.11) and its analog at real places, 
we have that 
$$
(w - a)^{j} \cdot  
\Cal{K}_{\infty}^{(j)}({\scriptstyle \frac{1}{2}} + it,  0, a,  \chi) \, \ge \, 0 \qquad (\text{for $w \le a$})
$$ 
Having all terms non-negative in (5.13) when $w < a$, we can interchange the first sum with 
the second and the integral. Since 
$$\Cal{K}_{\infty}({\scriptstyle \frac{1}{2}} + it,  0, w,  \chi) \, =\,
\sum_{j = 0}^{\infty}\,\, {(w - a)^{j}\over j!}  \,
\Cal{K}_{\infty}^{(j)}({\scriptstyle \frac{1}{2}} + it,  0, a,  \chi)
$$ 
the absolute convergence of (5.10) for $\Re(w) > 1$ follows.

Setting $w = 1 + \varepsilon$, then for arbitrary $T > 1,$ 
$$
\sum_{\chi \in \Chat_{0, S}}\;\;\;\,\int\limits_{\frak{I}_{_{\chi}}(T)}\,  
|L({\scriptstyle \frac{1}{2}} + it, f \otimes \chi)|^{2} \cdot T^{-1 - \varepsilon} \,
dt \,<\, Z(1 + \varepsilon) \, \ll_{\varepsilon}\, 1
$$ 
where $\frak{I}_{_{\chi}}(T)\,=\, \{t\in \Bbb{R}: \kappa_{_{\chi}}(t) \le T \}$, and hence 
$$
\sum_{\chi \in \Chat_{0, S}}\;\;\;\,\int\limits_{\frak{I}_{_{\chi}}(T)}\,  
|L({\scriptstyle \frac{1}{2}} + it, f \otimes \chi)|^{2}\,
dt \, \ll_{\varepsilon}\, T^{1 + \varepsilon} \tag 5.14
$$ 
Only finitely many characters contribute to the left-hand sum. This estimate is {\it compatible} 
with the convexity bound, in the sense that it implies for example that
$$
\int\limits_{0}^{T}\,  
|L({\scriptstyle \frac{1}{2}} + it, f)|^{2}\,
dt \, \ll_{\varepsilon}\, T^{^{[k: \Bbb{Q}] + \varepsilon}}
$$ 
Therefore, the function $Z(w)$ defined by (5.10) leads to averages of {\it reasonable} size 
suitable for applications. We return to a further study of the analytic properties of this function in a forthcoming paper.

\vskip20pt 
\noindent {\bf Concluding remarks:} The specific choice (3.11) of the data $\ph_{\nu}$ 
at archimedean places was made for {\it no} reason other than simplicity, enabling us 
to illustrate the non-vacuousness of the structural framework. Specifically, this choice 
allowed us to show that asymptotic formulas can be obtained, and that the averaging is 
{\it not too long}, i.e., compatible with the convexity bound. This choice sufficed for 
our purposes, which, again, were to stress generality, leaving aside the more 
technical aspects necessary in obtaining sharper results. 
Its use allowed us to quickly understand the size of the averages via the pole at $w = 1$, 
and dispensed with unnecessary complications.

The function $I(v, w)$ in (5.1) is analytic for $v$ in a neighborhood of $0$ and $\Re(w)$ 
sufficiently large. This follows easily from the analytic properties of $\Cal{K}_{\infty}(s, v, w, \chi)$ 
discussed in Section 5. By computing $I(v, w)$ using (4.13), this simple observation 
can be used to find the value of the constant $\kappa$ given in (4.7).



\vskip20pt\noindent {\bf \S Appendix 1. 
Convergence of Poincar\'e series
}
\vskip 10pt

\def\grad{\nabla}
\def\wt{g}

\def\less{\ll}

The aim of this appendix is to discuss the proofs of Proposition 2.6
and Theorem 2.7. Given the lack of complete arguments in the
literature, we have given a full account here, applicable more
generally. For a careful discussion of some aspects of $GL(2)$, see
\cite{GJ} and \cite{CPS1}.  Note that the latter source needs some
small corrections in the inequalities on pages 28 and 29.

We first prove the absolute convergence of the Poincar\'e series,
uniformly on compacts on $G_\A$, for $G=GL_2$ over a number field $k$
with ring of integers $\o$, for $\Re(\sv)>1$ and $\Re(\sw)>1$. Second,
we recall the notion of {\it norm} on a group, to prove convergence in
$L^2$ for {\it admissible} data (see the end of Section 2), also
reproving pointwise convergence by a more broadly applicable method.

Toward our first goal, we need an elementary comparison of sums and
integrals under mild hypotheses. Let $V_1,\ldots,V_n$ be
finite-dimensional real vector spaces, with fixed inner products, and
put  
$$
V \;=\; V_1 \oplus \ldots \oplus V_n
\disp{orthogonal direct sum}
$$
with the natural inner product. Fix a lattice $\Lambda$ in $V$, and
let $F$ be a period parallelogram for $\Lambda$ in $V$, containing
$0$. Let $g$ be a real-valued function on $V$ with 
$\wt(\xi)\ge 1$, such that $1/g$ has finite {\it integral} over
$V$, and is {\it multiplicatively bounded} on each
translate $\xi+F$, in the sense that, for each $\xi\in \Lambda$,
$$
\sup_{y\in \xi+F} {1\over \wt(y)} \;\less\; \inf_{y\in \xi+F}
{1\over \wt(y)} 
\disp{with implied constant independent of $\xi$}
$$ 
For a differentiable function $f$, let $\grad_i f$ be the gradient of $f$ in the $V_i$ variable. Then,
$$
\sum_{\xi\in \Lambda}\,\, |f(\xi)|
\;\less \;
\int_{V}\,\, |f(\xi)|\,d\xi
\; + \;
\sum_i 
\,\sup_{\xi\in V} 
\Big(
\wt(\xi)
\cdot |\grad_i f(\xi)|\Big)
$$
with the implied constant independent of $f$.

The following calculus argument gives this comparison (Abel summation).
Let ${\text{vol}}(\Lambda)$ be the natural measure of $V/\Lambda$. 
Certainly, 
$$
{\text{vol}}(\Lambda)\cdot \sum_{\xi\in \Lambda}\, |f(\xi)|
\;=\;
\sum_{\xi\in \Lambda}\,  |f(\xi)|\cdot \int_{\xi+F}\,dx
$$
and
$$
f(\xi)\int_{\xi+F} dx
\;=\;
\int_{\xi+F}(f(\xi)-f(x))\,dx
\;+\;
\int_{\xi+F}f(x)\,dx
$$
The sum over $\xi\in \Lambda$ of the latter integrals is obviously the
integral of $f$ over $V$, as in the claim. The differences
$f(\xi)-f(x)$ require further work. For $i=1,\ldots,n$, let $x_i$ and
$y_i$ be the $V_i$--components of $x,y\in V$, respectively. Let
$$
d_i(F) \,= \sup_{x,y\in F} |x_i - y_i|
$$
By the Mean Value Theorem, we have the easy estimate
$$
|f(\xi)-f(x)|
\,\le\,
\sum_{i=1}^n d_i(F) \cdot \sup_{y\in \xi+F} 
|\grad_i f(y)|
$$
Then,

\vskip5pt
\vbox{
$$
\sum_{\xi\in \Lambda}\, \,\int_{\xi+F}\,
|f(\xi)-f(x)|\,dx
\;\less\;
\sum_{\xi\in \Lambda}\;
\sum_{i=1}^n \;\sup_{y\in \xi+F} |\grad_i f(y)|
$$
$$
=\;
\sum_{i=1}^n \; \sum_{\xi\in \Lambda}\;
 \sup_{y\in \xi+F} 
\left({1\over \wt(y)}\, \wt(y) |\grad_i f(y)|
\right)
\;\le\;
\sum_{i=1}^n \; \sum_{\xi\in \Lambda}\,
\Big(\sup_{y\in \xi+F} {1\over \wt(y)} \Big)\cdot
\left(\sup_{y\in V} \wt(y)\,|\grad_i f(y)|
\right)
$$
$$
\less\;
\int_{V}\, {du\over \wt(u)}
\cdot \sum_i
\,\sup_{y\in V} \big(\wt(y) |\grad_i f(y)|\big)
\;\less\;
\sum_i
\,\sup_{y\in V} \big(\wt(y) |\grad_i f(y)|\big)
$$
}
\noindent This gives the indicated estimate.


The above estimate will show that the Poincar\'e series with parameter
$\sv$ is dominated by the sum of an Eisenstein series at $\sv$ and an
Eisenstein series at $\sv+1+\eps$ for every $\eps>0$, under mild
assumptions on the archimedean data. Such an Eisenstein series
converges absolutely and uniformly on compacts for $\Re(\sv)>1$,
either by Godement's criterion, in classical guise in \cite{B}, or by
more elementary estimates that suffice for $GL_2$. Thus, the
Poincar\'e series converges absolutely and uniformly for $\Re(\sv)>1$.

The assumptions on the archimedean data
$$
\Phi_\infty(x) \,=\, \ph_\infty\bmatrix{1&x\cr &1}
$$
are that
$$
\int_{k_\infty} |\Phi_\infty(\xi)|\,d\xi \,<\, +\infty
$$
and, letting $\grad_\v $ be the gradient along the summand $k_\v $ of
$k_\infty$, that, for some $\eps>0$, for each $\v|\infty$,
$$
\sup_{\xi\in k_\infty} (1+|\grad_\v  \Phi_\infty(\xi)|) \,<\, \infty
$$
The comparison argument is as follows.
To make a vector from which to form an Eisenstein series, left-average
the kernel 
$$
\ph\left(\bmatrix{a&\cr &d}\bmatrix{1 & x\cr &1}\right)
\,=\,
|a/d|^\sv \cdot \Phi(x)
\disp{extended by right $K_\A$--invariance}
$$
for the Poincar\'e series over $N_k$. That is, form
$$
\pht(g) \;=\; \sum_{\beta\in N_k}\ph(\beta\cdot g)
$$
This must be proven to be dominated by a vector (or vectors) from
which Eisenstein series are formed. The usual vector for standard
spherical Eisenstein series is 
$$
\eta_s\bmatrix{a& *\cr &d} \,=\, |a/d|^s
$$
extended to $G_\A$ by right $K_\A$--invariance. We claim that
$$
\pht \;\less\; \eta_{\sv} + \eta_{\sv+1+\eps}
\disp{for all $\eps>0$}
$$
Since all functions $\ph$, $\pht$ and $\eta_s$ are right
$K_\A$--invariant and have trivial central character, it suffices to
consider $g=nh$ with $n\in N_\A$ and 
$$
h \,=\, \bmatrix{y& \cr &1} \,\in\, H_\A
$$
Let
$$
n_t \,=\, \bmatrix{1 & t \cr &1}
$$
We have
$$
\ph(n_\xi\cdot n_x h)
\,=\,
\ph(h\cdot h^{-1}n_\xi n_xh)
\,=\,
\ph(h\cdot h^{-1}n_{\xi+x}h)
\,=\,
|y|^\sv \cdot \Phi\left({\textstyle{1\over y}}\cdot (\xi+x)\right)
$$
Thus, to dominate the Poincar\'e series by Eisenstein series, it
suffices to prove that 
$$
\sum_{\xi\in k} \, \Phi\left({\textstyle{1\over y}}\cdot (\xi+x)\right)
\,\less\,
1 \;+\; |y|
\disp{uniformly in $x\in N_\A$, $y\in \J$}
$$
Since $\pht$ is left $N_k$--invariant, it suffices to take $x\in\A$
to lie in a set of representatives $X$ for $\A/k$, such as
$$
X \;=\; k_\infty/\o \;\oplus\; {\textstyle\prod}_{\v<\infty}\, \o_\v 
$$
where, by abuse of notation, $k_\infty/\o$ refers to a period
parallelogram for the lattice $\o$ in $k_\infty$. As $\ph$ and $\pht$
are left $H_k$--invariant, so we can adjust $y$ in $\J$ by
$k^\times$. Since $\J^1/k^\times$ is compact, we can choose
representatives in $\J$ for $\J/k^\times$ lying in $C'\cdot (0,+\infty)$ for
some compact set $C'\subset \J^1$, with $(0,+\infty)$ embedded in $\J$
as usual by
$$
t \to (t^{1/n},t^{1/n},\ldots,t^{1/n},1,1,\ldots,1,\ldots)
\disp{non-trivial entries at archimedean places}
$$
where $n=[k:\Q]$. Further, for simplicity, we may adjust the
representatives $y$ such that $|y|_\v \le 1$ for all finite primes $\v$.
The compactness of $C'$ implies that
$$
|y|
\;\le\;
{\textstyle\prod}_{\v|\infty}\, |y_\v |_\v  
\;\less\;
|y|
\disp{with implied constant depending only on $k$}
$$
Likewise, due to the compactness, the archimedean valuations of
representatives have bounded ratios. 

At a finite place, $\Phi_\v \left({1\over y}\cdot (x+\xi)\right)$ vanishes unless
$$
{1\over y}\cdot (x+\xi) \,\in\, \o_\v 
$$
That is, since we want a uniform bound in $x\in\o_\v $, this vanishes
unless  
$$
\xi \,\in\, \o_\v  + y\cdot \o_\v  \subset \o_\v 
$$
since we have taken representatives $y$ with $y_\v $ integral at
all finite $\v$. Thus, the sum over $\xi$ in $k$ reduces to a sum over
$\xi$ with archimedean part in the lattice $\Lambda=\o\subset
k_\infty$.

Setting up a comparison as above, let 
$$
V\,=\, k_\infty\,=\;\bigoplus_{\v|\infty} k_\v 
$$
Let ${\text{vol}}(\Lambda)$ be the
volume of $\Lambda$. For archimedean place $\v$ let $\grad_\v $ be
the gradient along $k_\v $, and $d_\v (\Lambda)$ the maximum of
$|x_\v -y_\v |_\v $ for $x,y\in F$, a fixed period parallelogram for
$\Lambda$ in $k_\infty$. We have 
$$
\sum_{\xi\in \Lambda}\, 
\Phi\left({\textstyle{1\over y}}\cdot (\xi+x)\right)
\;\less \;
\int_{k_\infty} 
\Phi_\infty\left({\textstyle{1\over y}}\cdot (\xi+x)\right)
\, d\xi
\;+\;
\sum_{\v|\infty}\,\,\sup_{\xi\in k_\infty}
\Big(
\wt(\xi)
\cdot |\grad_\v  \Phi_\infty(\xi)|\Big)
$$
for any suitable weight function $\wt$. In the integral, replace $\xi$
by $\xi-x$, and then by $\xi\cdot y$, to see that
$$
\int_{k_\infty} 
\Phi_\infty\left({\textstyle{1\over y}}\cdot (\xi+x)\right)
 d\xi
\;=\;
|y|_\infty \cdot \int_{k_\infty} \Phi_\infty(\xi)\, d\xi
\;\less\;
|y| \cdot \int_{k_\infty} \Phi_\infty(\xi)\, d\xi
$$
with the implied constant depending only upon $k$, using the choice of
representatives $y$ for $\J/k^\times$.

To estimate the sum, for $x\in k_\infty$, fix $\eps>0$ and take weight
function  
$$
\wt(\xi) \,=\, {\textstyle\prod}_{\v|\infty}\, (1+|\xi|_\v ^2)^{\hf+\eps}
$$
This is readily checked to have the {\it multiplicative boundedness}
property needed: the function $\wt$ is continuous, and for $|\xi|\ge
2|x|$, we have the elementary
$$
\Hf\cdot |\xi| \,\le\, |\xi-x| \,\le\, 2\cdot |\xi|
$$
from which readily follows the bound for $\wt(\xi)$.

What remains is to compute the indicated supremums with attention to
their dependence on $y$. At an archimedean place $\v$, 

\vskip5pt
\vbox{
$$
\sup_{\xi\in k_\infty} \left(
\wt(\xi)
\cdot 
|\grad_i 
\Phi_\v \left({\textstyle{1\over y}}\cdot (\xi+x)\right)
|
\right)
\,=\,
\sup_{\xi\in k_\infty} \left(
\wt(\xi-x)
\cdot 
|\grad_i \Phi_\v \left({\textstyle{1\over y}}\cdot \xi\right)|
\right)
$$
$$
\;\less\;
\sup_{\xi\in k_\infty} \left(
\wt(\xi)
\cdot 
|\grad_i \Phi_\v \left({\textstyle{1\over y}}\cdot \xi\right)|
\right)
$$
}

\noindent by using the boundedness property of $\wt$. Then replace $\xi$ by
$\xi\cdot y$, to obtain
$$
\sup_{\xi\in k_\infty} \left(
\wt(y\cdot\xi)
\cdot 
|\grad_i \Phi_\v (\xi)|
\right)
$$
Since 
$$
(1+|y|^2_\v |\xi|^2_\v )\,\le\, 
(1+|y|^2_\v )\cdot (1+|\xi|^2_\v )
\disp{for all $\v|\infty$}
$$
we have $\wt(y\cdot \xi)\le \wt(y)\cdot \wt(\xi)$, and
$$
\sup_{\xi\in k_\infty} \left(
\wt(y\cdot \xi)
\cdot 
|\grad_\v  \Phi_\v (\xi)|
\right)
\,\,\le\,\,
\wt(y)\cdot \sup_{\xi\in k_\infty} \big(
\wt(\xi)
\cdot 
|\grad_i \Phi_\v (\xi)|
\big)
$$
Here the weighted supremums of the gradients appear, which we have
assumed finite.

Finally, estimate
$$
\wt(y)\,=\,\prod_{\v|\infty}\, (1+|y|_\v ^2)
$$
with $y$ in our specially chosen set of representatives. For these
representatives, for any two archimedean places $\v_1$ and $\v_2$, we
have 
$$
|y|^{n_{\v_{_1}}}_{\v_{_1}} \,\less\, |y|^{n_{\v_{_2}}}_{\v_{_2}}
$$
where the $n_{\v_i}$ are the local degrees $n_{\v_i}=[k_{\v_i}:\R]$. Therefore,
$$
|y|_\v  \,\less\, |y|^{n_\v /n}
$$
where $n=\sum_\v  n_\v $ is the global degree. Thus, 
$$
\prod_{\v|\infty}\, (1+|y|_\v ^2)
\,\less\,
1+|y|^2
$$
Then,
$$
\prod_{\v|\infty}\, (1+|y|_\v ^2)^{\hf+\eps}
\,\less\,
(1+|y|^2)^{\hf+\eps}
$$

Putting this all together, for every $\eps>0$
$$
\pht\bmatrix{y& * \cr 0&1} 
\,\less\, 
|y|^\sv \cdot (1 + |y|^2)^{\hf+\eps}
\,=\,
|y|^\sv + |y|^{\sv+1+2\eps}
\,=\,
\eta_{\sv}\bmatrix{y& * \cr 0&1} 
 + \eta_{\sv+1+2\eps}\bmatrix{y& * \cr 0&1} 
$$
which is the desired domination of the Poincar\'e series by a sum of
Eisenstein series.

For the particular choice of archimedean data
$$
\Phi_\infty(\xi)\,=\,\prod_{\v|\infty} \, {1\over (1+|\xi|^2_\v )^{\sw/2}}
$$
the {\it integrability} condition is met when $\Re(\sw)>1$. Similarly,
the weighted supremums of gradients are {\it finite} for $\Re(\sw)>1$.

Altogether, this particular Poincar\'e series is absolutely convergent
for $\Re(\sv)>1+2\eps$ and $\Re(\sw)>1+\eps$, for every $\eps>0$. This
proves Proposition 2.6.


\def\g{{\fr g}}
\def\Ad{\hbox{Ad}}

\def\card{\hbox{card}}
\def\op{{\scriptstyle {\text {o\hskip-1ptp}}}}
\def\less{\ll}
\def\bull{\hfil\break{$\bullet\;$}}

\def\psiobar{\overline{\psi}_o}

\vskip10pt
\noindent {\it Soft convergence estimates on Poincar\'e series}: 
Now we give a different approach to convergence, more convenient for
proving square integrability of Poincar\'e series. It is more robust,
and does also reprove pointwise convergence, but gives a weaker result
than the previous more explicit approach. Let $G$ be a (locally
compact, Hausdorff, separable) unimodular topological group. Fix a
compact subgroup $K$ of $G$. A {\it norm} $g\to \|g\|$ on $G$ is a
positive real-valued {\it continuous} function on $G$ with properties

\vskip 5pt 
\indent \vbox{
\bull $\|g\|\ge 1$ and $\|g^{-1}\|=\|g\|$
\bull {\it Submultiplicativity}: $\|gh\|\le\|g\|\cdot\|h\|$
\bull $K$--{\it invariance}: for $g\in G$, $k\in K$,
$ \|k\cdot g\| = \|g\cdot k\| = \|g\| $
\bull {\it Integrability}: for sufficiently large $\sigma>0$, 
$$
\int_G \, \|g\|^{-\sigma}\,dt \,<\, +\infty
$$
}

For a {\it discrete} subgroup $\Gamma$ of $G$, for $\sigma>0$ large enough
such that $\|g\|^{-\sigma}$ is integrable on $G$, we claim the
corresponding {\it summability}:
$$
\sum_{\gam\in\Gamma}{1\over \|\gamma\|^\sigma} \,<\, +\infty
$$
The proof is as follows. From
$$
\|\gam \cdot g\| \,\le\, \|\gam\|\cdot \|g\|
$$
for $\sigma>0$ 
$$
{1\over \|\gam\|^\sigma\cdot \|g\|^\sigma}
 \,\le\,
{1\over \|\gam \cdot g\|^\sigma}
$$
Invoking the discreteness of $\Gamma$ in $G$, let $C$ be a small
open neighborhood of $1\in G$ such that
$$
C \cap \Gamma \,=\, \{1\}
$$
Then,
$$
\int_C \, {dg\over \|g\|^\sigma} \cdot \sum_{\gam\in \Gamma}\, {1\over \|\gamma\|^\sigma}
\;\le\;
\int_C \,\sum_{\gam\in \Gamma}\, {1\over \|\gam\cdot g\|^\sigma}
\,dg
\;=\;
\sum_{\gam\in\Gamma}\, \int_{\gam^{-1}C} \,{1\over \|g\|^\sigma}\,dg
\;\le\;
\int_G \,{dg \over \|g\|^\sigma}\;<\; +\infty
$$
This gives the indicated summability.
Let $H$ be a closed subgroup of $G$, and define a {\it relative} norm
$$
\|g\|_{_H}  \, = \inf_{h\in H\cap\Gamma} \|h\cdot g\|
$$
From the definition, there is the left $H\cap\Gamma$--invariance
$$
\|h\cdot g\|_{_H} =\, \|g\|_{_H}
\disp{for all $h\in H\cap \Gamma$}
$$
Note that $\|\;\|_{_H}$ depends upon the discrete subgroup $\Gamma$. 


\vskip10pt
\noindent {\it Moderate increase, sufficient decay}: 
Let $H$ be a closed subgroup of $G$. A left $H\cap \Gamma$--invariant
complex-valued function $f$ on $G$ is of {\it moderate growth modulo}
$H\cap \Gamma$, when, for sufficiently large $\sigma>0$,
$$
|f(g)| \,\less\, \|g\|_{_H}^\sigma
$$
The function $f$ is {\it rapidly decreasing modulo}
$H\cap \Gamma$ if  
$$
|f(g)| \,\less\, \|g\|_{_H}^{-\sigma}
\disp{for {\it all} $\sigma>0$}
$$
The function $f$ is {\it sufficiently rapidly decreasing modulo}
$H\cap \Gamma$ (for a given purpose) if
$$
|f(g)| \, \less\, \|g\|_{_H}^{-\sigma}
\disp{for {\it some} sufficiently large $\sigma>0$}
$$
Since $\|g\|_{_H}$ is an infimum, for $\sigma>0$ the power
$\|g\|^{-\sigma}$ is a supremum
$$
{1\over \|g\|^\sigma} \;= \sup_{h\in H\cap \Gamma} {1\over
\|hg\|^\sigma} 
$$

\vskip10pt
\noindent{\it Pointwise convergence of Poincar\'e series}: 
We claim that, for $f$ left $H\cap\Gamma$--invariant and sufficiently
rapidly decreasing mod $H\cap \Gamma$, the {\it Poincar\'e series}
$$
P_f(g) \;\;\,= \sum_{\gam\in (H\cap \Gamma)\\\Gamma} f(\gam\cdot g)
$$
{\it converges absolutely and uniformly} on compacts.
To see this, first note that, for all $h\in H\cap \Gamma$,
$$
\|\gam \|_{_H} 
\,\le\,
\|h\cdot \gam \|
\,=\,
\|h\cdot \gam g \cdot g^{-1}\|
\,\le\,
\|h\gam g\| \cdot \|g^{-1}\|
$$
Thus, taking the inf over $h\in H\cap \Gamma$,
$$
{\|\gam\|_{_H}\over \|g^{-1}\|}
\,\le\,
\|\gam \cdot g\|_{_H} 
$$
Thus, for $\sigma>0$,
$$
{1\over \|\gam\cdot g\|_{_H}^\sigma} 
\,\le\,
{\|g\|^\sigma \over \|\gam \|_{_H}^\sigma}
$$
and
\vskip5pt
\vbox{
$$
P_f(g) 
\;\;\,= 
\sum_{\gam\in (H\cap \Gamma)\\\Gamma} \, f(\gam\cdot g)
\;\;\less 
\sum_{\gam\in (H\cap \Gamma)\\\Gamma}\, {1\over \|\gam\cdot g\|_{_H}^\sigma }
$$
$$
\;\le\;
\|g\|^\sigma \cdot\!\!\!\sum_{\gam\in (H\cap \Gamma)\\\Gamma} \, {1\over \|\gam\|_{_H}^\sigma}
\;\le\;
\|g\|^\sigma \cdot\!\!\!\sum_{\gam\in (H\cap \Gamma)\\\Gamma} 
\;\sum_{h\in H\cap \Gamma}\, {1\over \|h\cdot\gam\|^\sigma}
\;=\;
\|g\|^\sigma \cdot\,\,\,\,\!\!\!\sum_{\gam\in \Gamma}\, 
{1\over \|\gam\|^\sigma}
\;\less\; \|g\|^\sigma 
$$
}
\noindent estimating a sup of positive terms by the sum, for
$\sigma>0$ sufficiently large so that the sum over $\Gamma$ converges. 

\vskip10pt
\noindent {\it Moderate growth of Poincar\'e series}: 
Next, we claim that Poincar\'e series are of moderate growth modulo
$\Gamma$, namely, that 
$$
P_f(g) \, \less\, \|g\|_{_\Gamma}^\sigma
\disp{for sufficiently large $\sigma>0$}
$$
Indeed, the previous estimate is uniform in $g$, and the left-hand
side is $\Gamma$--invariant. That is, for all $\gam\in \Gamma$, 
$$
P_f(g) 
\,=\,
P_f(\gam\cdot g)
\,\less\, 
\|\gam\cdot g\|^\sigma 
\disp{with implied constant independent of $g,\gam$}
$$
Taking the inf over $\gam$ gives the assertion.

\vskip10pt

\noindent{\it Square integrability of Poincar\'e series}:
Next, we claim that for $f$ left $H\cap \Gamma$--invariant and
sufficiently rapidly decreasing mod $H\cap \Gamma$, $P_f$ is {\it
square-integrable} on $\Gamma\\G$. Unwind, and use the assumed estimate
on $f$ along with the above-proven moderate growth of the Poincar\'e
series: 
$$
\int_{\Gamma\\G} |P_f|^2
\;=\;
\int_{(H\cap \Gamma)\\G} |f|\cdot |P_f|
\;\less\;
\int_{(H\cap \Gamma)\\G} \|g\|^{-2\sigma}_{_H} \cdot \|g\|_{_H}^\sigma\,dg
$$
Estimating a sup by a sum, and unwinding further,
$$
\int_{(H\cap \Gamma)\\G} \|g\|^{-\sigma}_{_H} \,dg
\;\le\;
\int_{(H\cap \Gamma)\\G} 
\;
\sum_{h\in (H\cap \Gamma)\\\Gamma}
\|h\cdot g\|^{-\sigma}\,dg
\;=\;
\int_G
\|g\|^{-\sigma}\,dg
\,<\, +\infty
$$
for large enough $\sigma>0$. This proves the square integrability of
the Poincar\'e series.

\vskip10pt
\noindent{\it Construction of a norm on $PGL_2(\A)$}:
We want a norm on $G=PGL_2(\A)$ over a number field $k$ that meets the
conditions above, including the integrability, with $K$ the image in
$PGL_2(\A)$ of the maximal compact  
$$
\prod_{\v\isom \R} \;O_2(\R)
\times
\prod_{\v\isom \C} \;U(2)
\times
\prod_{\v<\infty} \;GL_2(\o_\v )
$$
of $GL_2(\A)$. We take $\Gamma$ to be the image in $PGL_2(\A)$ of
$GL_2(k)$.   
Let $\g$ be the algebraic Lie algebra of $GL_2$ over $k$, so that,
at each place $\v$ of $k$, 
$$
\g_\v  \,=\, \{\hbox{$2$-by-$2$ matrices with entries in $k_\v $}\}
$$
Let $\rho$ denote the Adjoint representation of $GL_2$ on $\g$,
namely,
$$
\rho(g)(x) \,=\, gxg^{-1}
\disp{for $g\in GL_2$ and $x\in \g$}
$$
The kernel of $\rho$ on $GL_2$ is the center $Z$, so the image $G$ of 
$GL_2$ under $\rho$ is $PGL_2$. As expected, let
$$
G_\v  \,=\, \rho(GL_2(k_\v )) \,=\, GL_2(k_\v )/Z_\v 
\hskip30pt
K_\v  \,=\, \rho(GL_2(o_\v )) 
\,=\, GL_2(\o_\v )/(Z_\v \cap GL_2(\o_\v ))
$$
and
$$
\Gamma \,=\,G_k \,=\, \rho(GL_2(k))\,=\,GL_2(k)/Z_k
$$
Since $\Gamma$ is a subgroup of $GL_k(\g_k)$, it is discrete in the
adelization of $GL_k(\g_k)$, so is discrete in $G_\A$. Let
$\{e_{ij}\}$ be the $2$-by-$2$ matrices with non-zero entry just at
the $(i,j)^{\text{th}}$ location, where the entry is $1$.


At an archimedean place $\v$ of $k$, put a Hilbert space structure on
$\g_\v $ by   
$$
\bra x,y\ket 
\,=\,
\tr (y^*x)
$$
where $y^*$ is $y$--transpose for $\v$ real, and $y$--transpose-conjugate
for $\v$ complex. We put the usual (sup-norm) operator norm on
linear operators $T$ on $\g_\v$, namely
$$
|T|_\op
\,=\,
\sup_{|x|\le 1} |Tx|
$$
By design, since the inner product on $\g_\v $ is $\rho(K_\v )$--invariant,
this operator norm is invariant under $\rho(K_\v )$.  


For a non-archimedean local field $k$ with norm $|\,\cdot \,|_\v $ and ring of
integers $\o$, give $\g_\v $ the sup-norm 
$$
|{\textstyle\sum}_{ij} a_{ij}\,e_{ij}|\,=\, \sup_{ij} |a_{ij}|_\v 
\disp{with $a_{ij}\in k_\v $}
$$
There is the {\it operator norm} on $GL_{k_\v }(\g_\v )$ given by 
$$
|g|_\op
\;\;= 
\sup_{x\in V,\;|x|\le 1} |g\cdot x|
$$
By design, this norm is invariant under $\rho(K_\v )$.


\vskip10pt
\noindent{\it Norms on local groups and adele groups}:
For any place $\v$ of $k$, define a (local) {\it norm} $\|g\|_\v $ on the
image $G_\v =PGL_2(k_\v )$ of $GL_2(k_\v )$ in $GL_{k_\v }(\g_\v )$ by 
$$
\|g\|_\v  \,=\, \max\{ |g|_\op,\;|g^{-1}|_\op\}
$$
Since the norm on $\g_\v $ is $K_\v $--invariant, and $K_\v $ is stable under
inverse, the operator norms are left and right
$K_\v $--invariant, and the norms $\|\;\|_\nu$ are left and right
$K_\v $--invariant. 
Note that for $\v<\infty$ the operator norm is $1$ on $K_\v $. To prove
that
\vskip3pt
\vbox{
$$
\|g\cdot h\|_\v  \, \le\, \|g\|_\v \cdot \|h\|_\v 
$$	
use the definition:
$$
\|g\cdot h\|_\v  
\,=\,
\max\{|gh|_\op,\;|h^{-1}g^{-1}|_\op\}
$$
$$
\,\le\,
\max\{|g|_\op\cdot |h|_\op,\;|g^{-1}|_\op\cdot |h^{-1}|_\op\}
\,\le\,
\max\{|g|_\op,\;|g^{-1}|_\op\}
\cdot
\max\{|h|_\op,\;|h^{-1}|_\op\}
\,=\,
\|g\|_\v \cdot \|h\|_\v 
$$
}
\noindent For $g=\{g_\v \}$ in the adele group $G_\A$, let
$$
\|g\| \,=\, {\textstyle\prod}_\v  \, \|g_\v \|_\v 
$$
The factors in the product are $1$ for all but finitely many places.
The left and right $K$--invariance for $K=\prod_\v  K_\v $ follows from the
local $K_\v $--invariance. Invariance under inverse is likewise clear.


\vskip10pt
\noindent{\it Integrability}:
Toward integrability, we explicitly bound the local integrals
$$
\int_{G_\v } \|g\|_\v ^{-\sigma} \, dg
$$
At finite primes, use the $p$--adic Cartan decomposition (here just
the elementary divisor theorem) inherited from $GL_2(k_\v )$ via the
quotient map, namely, 
$$
G_\v  \;\;\;\;= \bigsqcup_{\delta\in A_\v /(A_\v \cap Z_\v K_\v )} K_\v \cdot \delta \cdot K_\v 
\disp{where $A_\v $ is diagonal matrices}
$$
By conjugating by permutation matrices and adjusting by $Z_\v $, we may
assume, further, that 
$$
\delta\,=\, \bmatrix{\delta_1 & \cr & 1}
\disp{with $|\delta_1|\ge 1$}
$$
For any choice $\varpi_\v $ of local parameter for $k_\v $, we may adjust
by $A_\v \cap K_\v $ so that $\delta_1$ is a power of $\varpi$. On a given 
$K_\v $ double coset, the norm is
$$
\| K_\v \cdot \delta \cdot K_\v \|_\v 
\,=\,
\|\delta\|_\v 
\,=\, \max\{|\rho(\delta)|_\op,\; |\rho(\delta^{-1})|_\op\}
\,=\,
\max \{ |\delta_1|_\v,\;|\delta_1|^{-1}_\v  \}
$$
and
$$
\meas(K_\v \delta K_\v )
\,=\,
\meas(K_\v ) \cdot \card(K_\v \\K_\v \delta K_\v )
$$
Let $q=q_\v $ be the residue field cardinality, and let
$|\delta|_\v =q^\ell$ with $\ell\ge 0$. Then,        
$$
\card\;K_\v \\K_\v \delta K_\v 
\,=\,
\card\;(K_\v \cap \delta^{-1}K_\v \delta)\\ K_\v 
\,\le\,
\card\;K_\v (\ell)\\ K_\v 
$$
where $K_\v (\ell)$ is a sort of congruence subgroup, namely,
$$
K_\v (\ell) \,=\, \left\{ \bmatrix{a&b\cr c&d}\in K_\v \,:\, c \in \varpi^\ell\cdot \o_\v 
\right\}
$$
Let $K_\v' =\{g\in K_\v :g=I_2 \!\mod \varpi \o\}$, and let
$\F_q$ be the finite field with $q$ elements. We have an elementary
estimate
$$
[K_\v :K_\v (1)] 
\,=\, 
{[K_\v :K_\v ']\over [K_\v (1):K_\v ']}
\,=\,
\card\{ \hbox{lines in $\F_q^2$} \}
\,=\,
{q^2-1\over q-1}
\,=\,
q+1
\,<\,
q^2
$$
and
$$
[K_\v (\ell):K_\v (\ell+1)] \,=\, q
\disp{for $\ell\ge 1$}
$$
Thus, 
$$
[K_\v :K_\v (\ell)] \,\le\, q^2 \cdot q^{\ell-1}
\disp{for $\ell\ge 1$}
$$
Thus, the integral of $\|g\|_\v ^{-\sigma}$ has an upper bound
$$
\int_{G_\v } {dg\over \|g\|_\v ^\sigma}
\,\le\,
1 + \sum_{\ell\ge 1} q^{-\sigma\ell}\cdot q^{2 + (\ell-1)}
\,\le\,
1 + q\cdot \sum_{\ell\ge 1} (q^{1-\sigma})^\ell
$$
For $\sigma>1$, the geometric series converges. Thus,
$$
\int_{G_\v } {dg\over \|g\|_\v ^\sigma} 
\,\le \,
1 + q\cdot { q^{1-\sigma} \over 1 - q^{1-\sigma}}
\,<\,
{
1 + q\cdot q^{1-\sigma}
\over 
1 - q^{1-\sigma} 
}
\,=\,
{
1 + q^{2-\sigma}
\over 
1 - q^{1-\sigma} 
}
$$
Note that there is {\it no} leading constant. 


The integrability condition on the adele group can be verified by
showing the finiteness of the product of the corresponding local
integrals. Since there are only finitely many archimedean places, it
suffices to consider the product over finite places. By comparison to
the zeta function of the number field $k$, a product 
$$
{\underset{\nu < \infty}\to{\prod}}\,
{1+q_\v ^{-a} \over 1-q_\v ^{-b}}
\;=\,
{\underset{\nu < \infty}\to{\prod}}\,
{1-q_\v ^{-2a} \over (1-q_\v ^{-b})(1-q_\v ^{-a})}
\,=\,
{\zeta_k(a)\cdot \zeta_k(b)
\over
\zeta_k(2a)
}
$$
converges for $a>1$ and $b>1$. Thus, letting $G_{\text{fin}}$ be the
finite-prime part of the idele group $G_\A$, 
$$
\int_{G_{\text{fin}}} {dg\over \|g\|^\sigma} 
\;=\,
{\underset{\nu < \infty}\to{\prod}}\,
\int_{G_\v } {dg\over \|g\|_\v ^\sigma} 
\,<\,+\infty
$$
for $\sigma>0$ sufficiently large, from the previous estimate on the
corresponding local integrals.


For integrability locally at archimedean places, exploit the left and
right $K_\v $--invariance, via Weyl's integration 
formula. Let $A_\v $ be the image under $\Ad$ of the standard maximal
{\it split} torus from $GL_2(k_\v )$, namely, real diagonal matrices. Let
$\Phi^+=\{\alf\}$ be the singleton set of standard positive roots
of $A_\v $, namely 
$$
\alf\;:\;\bmatrix{a_1 & \cr & a_2} \to a_1/a_2
$$
$\g_\alf$ be the $\alf$--rootspace, and, for $a\in A_\v $,
let 
$$
D(a) \,=\, |\alf(a) - \alf^{-1}(a)|^{\dim_{_\R} \g_\alf}
$$
The Weyl formula for a left and right $K_\v $--invariant function $f$ on
$G_\v $ is 
$$
\int_{G_\v } f(g) \,dg
\,=\,
\int_{A_\v } D(a)\cdot f(a) \,da
$$
For $PGL_2$, the dimension $\dim_{_\R}\g_\alf$ is $1$ for $k_\v \isom \R$
and is $2$ for $k_\v \isom \C$. The norm of a diagonal element is easily
computed via the adjoint action on $\g_\v $, namely 
$$
\|a\|_\v  \,=\, \max\{|a_1/a_2|,\, |a_2/a_1|\}
$$
with the usual absolute value on $\R$. Thus, 
$$
D(a) \,\less\,\|a\|_\v ^{d_\v } 
\disp{with $d_\v =[k_\v :\R]$}
$$
Thus, the integral over $PGL_2(k_\v )$ is dominated by a one-dimensional
integral, namely, 
$$
\int_{G_\v } {dg\over \|g\|^\sigma_\v }
\,=\,
\int_{A_\v } {D(a)\over \|a\|_\v ^\sigma} \, da
\,\less\,
\int_{\R^\times} (\max(|x|,|x|^{-1})^{d_\v -\sigma}\,dx
\disp{with $d_\v =[k_\v :\R]$}
$$ 
The latter integral is evaluated in the fashion
$$
\int_{\R^\times} (\max(|x|,|x|^{-1})^{-\beta}\,dx
\,=\,
\int_{|x|\le 1} (|x|^{-1})^{-\beta}\, dx
+
\int_{|x|\ge 1} |x|^{-\beta}\, dx
\,<\,
+\infty
$$
for $\v$ either real or complex. This gives the desired local
integrability for large $\sigma$ at archimedean places, and completes
the proof of global integrability.

\vskip10pt
\noindent{\it Poincar\'e series for $GL_2$}: 
Recall the context of Sections 2 and 3. Let $G=GL_2(\A)$ over a number
field $k$, $Z$ the center of $GL_2$, and $K_\v $ the standard maximal
compact in $G_\v $. Let
$$
M\,=\, \left\{\bmatrix{*&0\cr 0&*}\right\}
\hskip20pt
N\,=\, \left\{\bmatrix{1&*\cr 0&1}\right\}
$$
To form a Poincar\'e series, let
$\ph=\bigotimes_\v \ph_\v $, where each $\ph_\v $ is right $K_\v $--invariant,
$Z_\v $--invariant, and on $G_\v $ 
$$
\ph_\v \left( \bmatrix{a& \cr &1}\bmatrix{1 & x \cr &1}\right)
\,=\,
|a|_\v ^\sv \cdot \Phi_\v (x)
$$
where at finite primes $\Phi_\v $ is the characteristic function of the
local integers $\o_\v $. At archimedean places, we assume that $\Phi_\v $
is sufficiently continuously differentiable, and that these
derivatives are absolutely integrable. The global function $\ph$ is
left $M_k$--invariant, by the product formula. Then, let
$$
f(g) \,=\, \int_{N_\A} \psibar(n)\, \ph(ng)\,dn
$$
where $\psi$ is a standard non-trivial character on $N_k\\N_\A\isom
k\\\A$. As in (4.6), but with slightly different notation, the
Poincar\'e series of interest is 
$$
\Pepsi(g) \;\;\;= \sum_{\gam\in Z_kN_k\\G_k} f(\gam\cdot g) 
$$

\vskip10pt
\noindent{\it Convergence uniformly pointwise and in $L^2$}:
From above, to show that this converges absolutely and uniformly on
compacts, and also that it is in $L^2(Z_\A G_k\\G_\A)$, use a norm on
the group $PGL_2=GL_2/Z$, take $\Gamma=PGL_2(k)$, and show that $f$ is
{\it sufficiently rapidly decreasing on $PGL_2(\A)$ modulo $N_k$.} 

To give the sufficient decay modulo $N_k$, it suffices to prove
sufficient decay of $f(nm)$ for $n$ in a well-chosen set of
representatives for $N_k\\N_\A$, and for $m$ in among representatives 
$$
m \,=\, \bmatrix{a & \cr &1}
$$
for $M_\A/Z_\A$. For $m\in M_\A$ and $n\in N_\A$, the
submultiplicativity $\|nm\|\le\|n\|\cdot \|m\|$ gives
$$
{1\over \|n\|^\sigma\cdot \|m\|^\sigma}
\,\le\,
{1\over \|nm\|^\sigma}
\disp{for $\sigma>0$}
$$
That is, roughly put, it suffices to prove decay in $N_\A$ and $M_\A$
separately. Since $N_k\\N_\A$ has a set of representatives $E$ that is 
{\it compact}, on such a set of representatives the norm is {\it
bounded}. Thus, it suffices to prove that
$$
f(nm) \,\less\, {1\over \|m\|^\sigma} \qquad \quad
\left(\text{for $n\in E$, and $m=\bmatrix{a& \cr &1}$}\right)
$$
Since $f$ factors over primes, as does $\|m\|$, it suffices to give
suitable {\it local} estimates.


At finite $\v$, the $\v^{\text{th}}$ local factor of $f$ is left
$\psi$--equivariant by $N_\v $, and  
\vskip5pt
\vbox{
$$
f_\v (nm)\,=\,
\psi(n)\cdot
\int_{N_\v } \psibar(n') \, \ph_\v (n'm)\,dn'
\,=\,
\psi(n)\cdot
\int_{N_\v } \psibar(n') \, \ph_\v (m\cdot m^{-1}n'm)\,dn'
$$
$$
\,=\,
\psi(n)\cdot
|a|_\v \cdot \int_{k_\v } \psiobar(ax) \, |a|_\v ^\sv \cdot \Phi_\v (x)\,dx
$$
}
\noindent where
$$
\psi\bmatrix{1 & x \cr & 1}\;=\; \psiobar(x)
\hskip30pt
m\,=\, \bmatrix{a & \cr & 1}
$$
Then,
$$
|f_\v (nm)| \,=\,
|a|_\v ^{\sv+1}\cdot \int_{k_\v } \psiobar(ax) \, \Phi_\v (x)\,dx
\,=\,
|a|_\v ^{\sv+1}\cdot \Phat_\v (a)
$$


At every finite place $\v$, $\Phi_\v $ has compact support, and at almost
every finite $\v$, $\Phat_\v $ is simply the characteristic function of
$\o_\v $. Thus, almost everywhere, 
$$
|f_\v (nm)| \,\le\,
|a|_\v ^{\sv+1}\cdot \Phat_\v (a)
\,\le\,
\left(\max\{|a|_\v ,\,|a|_\v ^{-1}\}\right)^{-(\sv+1)}
\,=\,
\|m\|_\v ^{-(\sv+1)}
$$
At the finitely many finite places where $\Phat$ is not exactly the
characteristic function of $\o_\v $, the same argument still gives the
weaker estimate
$$
|f_\v (nm)| \,\le\,
|a|_\v ^{\sv+1}\cdot \Phat_\v (a)
\,\less\,
\left(\max\{|a|_\v ,\;|a|_\v ^{-1}\}\right)^{-(\sv+1)}
\,=\,
\|m\|_\v ^{-(\sv+1)}
$$
Thus, we have the finite-prime estimate
$$
{\underset{\nu < \infty}\to{\textstyle\prod}}\, |f_\v (nm)|
\,\less\,
{\underset{\nu < \infty}\to{\textstyle\prod}}\,
\|m\|_\v ^{-(\sv+1)}
$$

At archimedean places, given $\ell>0$, for $\Phi_\v $ sufficiently
differentiable with absolutely integrable derivatives, ordinary
Fourier transform theory implies that
$$
|\Phat_\v (a)|\,\less\,(1+|a|_\v )^{-\ell}
$$
Thus, from the general local calculation above,
$$
|f_\v (nm)|
\,=\,
|a|_\v ^{\sv+1}\cdot \Phat_\v (a)
\,\less\,
|a|_\v ^{\sv+1}\cdot (1+|a|_\v )^{-\ell}
$$
This gives the sufficient decay of $f$ at archimedean places.

In summary, since the local factors $f_\v $ of $f$ have sufficient
decay, the function $f$ has sufficient decay so that the associated
Poincar\'e series $\Pepsi=P_f$ converges uniformly on compacts, and is
in $L^2(Z_\A GL_2(k)\\GL_2(\A))$. This proves Theorem 2.7.


\vskip20pt\noindent {\bf \S Appendix 2. Mellin transform of Eisenstein
Whittaker functions}
\vskip 10pt 

The computation discussed in this appendix was needed in the proof of Proposition 4.10. 
While the details of this computation are given below, we also cite
\cite{W2}, Chapter VII, for standard facts about the Tate-Iwasawa theory
of zeta integrals.

The global Mellin transform of $W^E$ factors
$$
\int_{\J} |a|^{v}\, W^E_{s,\,\chi}\pmatrix a&0\cr 0&1\endpmatrix da
\;=\;
\prod_{\nu} \int_{k_{\nu}^{\times}} |a|_{\nu}^{v}\, W^E_{s,\,\chi,\, \nu}\pmatrix a&0\cr 0&1\endpmatrix da
$$
To compute this, we cannot simply change the order of integration,
since this would produce a divergent integral along the way.
Instead, we present the vectors $\eta_{\nu}$ in a different form.
Let $\Phi_{\nu}$ be any Schwartz function on $k_{\nu}^2$ invariant under 
$K_{\nu}$ (under the obvious right action of $GL_2$), and put
$$
\eta_{\nu}'(g) 
= \chi_{\nu}(\det g)\,|\det g|_{\nu}^{s} \cdot 
\int_{k_{\nu}^{\times}} \chi_{\nu}^{2}(t) \,|t|_{\nu}^{2s} \cdot \Phi_{\nu}(t\cdot e_2 \cdot g)\,dt  
$$
where $e_{2} = e_{2, \, \nu}$ is the second basis element in $k_{\nu}^{2}.$ This $\eta_{\nu}'$ has
the same left $P_{\nu}$--equivariance as $\eta_{\nu}$, namely 
$$
\eta_{\nu}' \left(\pmatrix a&*\cr 0&d\endpmatrix\cdot g \right)
\,=\,
|a/d|_{\nu}^{s}\cdot \chi_{\nu}(a/d)\cdot 
\eta_{\nu}'(g)$$
For $\Phi_{\nu}$ invariant under the standard maximal compact $K_{\nu}$ of
$GL_{2}(k_{\nu})$, the function $\eta_{\nu}'$ is right $K_{\nu}$--invariant. By the
Iwasawa decomposition, up to constant multiples, there is only one such
function, so 
$$
\eta_{\nu}'(g) = \eta_{\nu}'(1) \cdot \eta_{\nu}(g)
\disp{since $\eta_{\nu}(1)=1$}
$$
and\footnote{From now on, to avoid clutter, suppress the subscript $\nu$ 
where there is no risk of confusion. For instance, we shall write $|\cdot|,$ $\psi,$ $\chi,$ etc., 
rather than $|\cdot|_{\nu},$ $\psi_{\nu},$ $\chi_{\nu},$ etc.}
$$
\eta_{\nu}'(1) 
= \int_{k_{\nu}^{\times}} \chi^{2}(t)\, |t|^{2s} \cdot \Phi(t\cdot e_{2} \cdot 1)\,dt
= \zeta_{\nu}(2s,\, \chi^2,\, \Phi(0,*))\qquad
(\text{a Tate-Iwasawa zeta integral})
$$
Thus, it suffices to compute the local Mellin transform of

\vskip5pt
\vbox{
$$
\eta'_{\nu}(1)\cdot W^{E}_{s,\, \chi,\, \nu}(m)
=
\int_{N_{\nu}} \psibar(n)\,\eta'_{\nu}(w_{\circ}  nm)\,dn
= 
\chi(a)|a|^s
\int_{N_{\nu}} \psibar(n)\,
\int_{k_{\nu}^{\times}} \chi^2(t)|t|^{2s}\,\Phi(t\cdot e_2\cdot w_{\circ}  nm)
\,dt\,dn
$$
$$
= 
\chi(a)|a|^s
\int_{k_{\nu}} \psibar(x)\,
\int_{k_{\nu}^\times} \chi^2(t)|t|^{2s}\,\Phi(tx, ta)
\,dt\,dx
\qquad \left(\text{with $m=\pmatrix a&0\cr 0&1\endpmatrix$}\right)
$$
}
At {\it finite} primes $\nu$, we may as well take $\Phi$ to be
$$
\Phi(t, x) = \ch_{\o_{\nu}}(t) \cdot \ch_{\o_{\nu}}(x)
\disp{$\ch_{X}=$ characteristic function of set $X$}
$$
Then $\eta_{\nu}'(1)$ is exactly an $L$--factor (see \cite{W2}, page 119, Proposition 10)
$$
\eta_{\nu}'(1) 
= \zeta_{\nu}(2s,\, \chi^2,\, \ch_{\o_{\nu}})
= L_{\nu}(2s,\, \chi^2)
$$
and (for further details see \cite{W2}, page 107, Corollary 1, and page 108, Corollary 3), 

\vskip5pt
\vbox{
$$
\eta'_{\nu}(1)\cdot W^{E}_{s,\,\chi,\, \nu}\pmatrix a&0\cr 0&1\endpmatrix
=
\chi(a)|a|^s
\int_{k_{\nu}} \psibar(x)\,\ch_{\o_{\nu}}(tx)
\int_{k_{\nu}^{\times}} \chi^2(t)|t|^{2s}\,\ch_{\o_{\nu}}(ta)
\,dt\,dx
$$
$$
=
\chi(a)|a|^s
\,\meas(\o_{\nu})\,
\int_{k_{\nu}^{\times}} \ch_{\o_{\nu}^*}(1/t) \chi^2(t)|t|^{2s-1}\,\ch_{\o_{\nu}}(ta)
\,dt
$$
$$
=
|\diff_{_\nu}|^{1/2}
\cdot
\chi(a)|a|^s
\int_{k_{\nu}^{\times}} \ch_{\o_{\nu}^*}(1/t) \chi^2(t)|t|^{2s-1}\,\ch_{\o_{\nu}}(ta)
\,dt
$$
}
\noindent where $\diff_{\nu}\in k_{\nu}^{\times}$ is such that $(\o_{\nu}^*)^{-1}=\diff_{\nu}\cdot
\o_{\nu}$. We can compute now the Mellin transform
$$
\int_{k_{\nu}^{\times}} |a|^{v} \cdot 
\left(
\chi(a)|a|^s
\int_{k_{\nu}^{\times}} \ch_{\o_{\nu}^*}(1/t) \chi^2(t)|t|^{2s-1}\,\ch_{\o_{\nu}}(ta)
\,dt\right)
\,da
$$
Replace $a$ by $a/t$, and then $t$ by $1/t$ to obtain a product
of two zeta integrals   

\vskip4pt
\vbox{
$$
\left(
\int_{k_{\nu}^{\times}} |a|^{v} \cdot \chi(a)|a|^s\,\ch_{\o_{\nu}}(a)\,da\right)
\cdot 
\left(
\int_{k_{\nu}^{\times}} \ch_{\o_{\nu}^*}(1/t) \chi(t)|t|^{s-1-v}\,
\,dt\right)
$$
$$
= 
\zeta_{\nu}(v+s,\, \chi,\, \ch_{\o_{\nu}})
\cdot
\zeta_{\nu}(v+1-s,\, \chibar,\, \ch_{\o_{\nu}^{*-1}})
$$
$$
= 
L_{\nu}(v+s,\, \chi)
\cdot 
L_{\nu}(v + 1 - s,\, \chibar)
\cdot 
|\diff_{_\nu}|^{-(v + 1 - s)}\, \chi(\diff_{_\nu})
$$
}
Thus, dividing through by $\eta'_{\nu}(1)$ and putting back the measure
constant, the Mellin transform of $W^{E}_{s,\, \chi,\, \nu}$ is 

\vskip5pt
\vbox{
$$
\int_{k_{\nu}^{\times}} |a|^{v}\,W^{E}_{s,\,\chi,\, \nu}\pmatrix a&0\cr 0&1\endpmatrix\,da
= 
|\diff_{_\nu}|^{1/2}\cdot
{
L_{\nu}(v+s,\, \chi)
\cdot 
L_{\nu}(v + 1 - s,\chibar)
\over
L_{\nu}(2s, \, \chi^2)
}
\cdot 
|\diff_{_\nu}|^{-(v + 1 - s)}\,\chi(\diff_{_\nu})
$$
}
Let $\diff$ be the idele whose $\nu^{\text{th}}$ component is $\diff_{\nu}$ for
finite $\nu$ and whose archimedean components are all $1$. The product
over all finite primes $\nu$ of these local factors is    
$$
\int_{\J^{\text{fin}}} |a|^{v}\,W^{E}_{s,\, \chi}\pmatrix a&0\cr 0&1 \endpmatrix\,da
=
|\diff|^{1/2}
\cdot
{L(v+s,\, \chi)\cdot L(v+1-s,\, \chibar)
\over
L(2s,\, \chi^2)
}
\cdot 
|\diff|^{-(v + 1 - s)}\,\chi(\diff)
$$
In our application, we will replace $s$ by $1 - s$ and $\chi$ by
$\chibar$, giving
$$
\int_{\J^{\text{fin}}} |a|^{v}\,W^{E}_{1-s,\, \chibar}\pmatrix a&0\cr 0&1 \endpmatrix\,da
=
|\diff|^{1/2}
\cdot
{L(v+1-s,\, \chibar)\cdot L(v+s,\, \chi)
\over
L(2-2s,\, \chibar^2)
}
\cdot 
|\diff|^{-(v + s)}\,\chibar(\diff)
$$
In particular, with $\chi$ trivial, 
$$
\int_{\J^{\text{fin}}} |a|^{v}\,W^{E}_{1-s}\pmatrix a&0\cr 0&1 \endpmatrix\,da
=
|\diff|^{1/2}
\cdot
{\zeta_{k}(v+1-s)\cdot \zeta_{k}(v+s)
\over
\zeta_{k}(2-2s)
}
\cdot 
|\diff|^{-(v+s)}
$$


\enddocument